\documentclass[a4paper,10pt]{article}
\usepackage{latexsym}
\usepackage{amsmath}
\usepackage{amssymb}
\usepackage{amscd}

\newlength{\minitwocolumn}
\makeatletter
\@addtoreset{equation}{section}
\makeatother

\title{Vertex operator approach for 
correlation functions of \\ 
Belavin's $(\mathbb{Z}/n\mathbb{Z})$-symmetric model}

\author{Yas-Hiro Quano\thanks
{email: quanoy@suzuka-u.ac.jp}}

\date{\it Department of Clinical Engineering, 
Suzuka University of Medical Science \\
      \it Kishioka-cho, Suzuka 510-0293, Japan \\
23 Oct 2008 \\[1cm]
Dedicated to Professor Tetsuji Miwa on the occasion of his sixtieth birthday}
\begin{document}

\maketitle
\begin{abstract}
Belavin's $(\mathbb{Z}/n\mathbb{Z})$-symmetric model is considered 
on the basis of bosonization of vertex operators 
in the $A^{(1)}_{n-1}$ model and vertex-face transformation. 
The corner transfer matrix (CTM) Hamiltonian of 
$(\mathbb{Z}/n\mathbb{Z})$-symmetric model and 
tail operators are expressed in terms of bosonized vertex 
operators in the $A^{(1)}_{n-1}$ model. Correlation functions 
of $(\mathbb{Z}/n\mathbb{Z})$-symmetric model can be 
obtained by using these objects, in principle. 
In particular, we calculate spontaneous polarization, 
which reproduces the result by ourselves in 1993.
\end{abstract}

\section{Introduction}

In this paper we consider 
Belavin's $(\mathbb{Z}/n\mathbb{Z})$-symmetric model \cite{Bela} 
on the basis of bosonization of vertex operators 
in the $A^{(1)}_{n-1}$ model \cite{AJMP} 
and vertex-face transformation. 
Belavin's $(\mathbb{Z}/n\mathbb{Z})$-symmetric model 
is a higher rank generalization of 
Baxter's eight-vertex model \cite{ESM} in the sense that 
the former model is an $n$-state model. The 
$(\mathbb{Z}/n\mathbb{Z})$-symmetric model is a vertex model 
on a two dimensional lattice such that the state variables 
take on values of $(\mathbb{Z}/n\mathbb{Z})$-spin. A local weight 
$R^{ik}_{jl}$ is assigned to spin configulation $j,l,i,k$ 
around a vertex. The model is $(\mathbb{Z}/n\mathbb{Z})$-symmetric 
in a sense that $R^{ik}_{jl}$ satisfies the two conditions: 
(i) $R^{ik}_{jl}=0$ unless $j+l=i+k$ (mod $n$), and (ii) 
$R^{i+p\,k+p}_{j+p\,l+p}=R^{ik}_{jl}$ for any $p\in (\mathbb{Z}/n\mathbb{Z})$. 
Since there are $n^3$ non zero weights among $R^{ik}_{jl}$'s, 
we may call the $(\mathbb{Z}/n\mathbb{Z})$-symmetric model by 
the $n^3$-vertex model. (When $n=2$, it becomes the eight-vertex model.) 

In \cite{LaP} Lashkevich and Pugai presented the integral formulae 
for correlation functions of the eight-vertex model \cite{ESM} using 
bosonization of vertex operators 
in the eight-vertex SOS model \cite{LuP} 
and vertex-face transformation. The present paper aims to give an 
$sl(n)$-generalization of Lashkevich-Pugai's construction. 
For our purpose we use the vertex-face correspondence between 
the $(\mathbb{Z}/n\mathbb{Z})$-symmetric model and 
unrestricted $A^{(1)}_{n-1}$ model. 
We first notice that the $A^{(1)}_{n-1}$ model \cite{JMO} 
is a restricted model, while we should relate 
the $(\mathbb{Z}/n\mathbb{Z})$-symmetric model with 
unrestricted $A^{(1)}_{n-1}$ model. We second notice that 
the original vertex-face correspondence \cite{JMO} maps 
the $A^{(1)}_{n-1}$ model in regime III to 
$(\mathbb{Z}/n\mathbb{Z})$-symmetric model in the 
disordered phase. We should relate the former with 
$(\mathbb{Z}/n\mathbb{Z})$-symmetric model in the 
principal regime. 

In this paper we present integral formulae for correlation 
functions of $(\mathbb{Z}/n\mathbb{Z})$-symmetric model 
on the basis of free field representation formalism. 
As the simplest example, we perform the calculation of 
the integral formulae for one-point function, in order to 
obtain the spontaneous polarization of 
$(\mathbb{Z}/n\mathbb{Z})$-symmetric model. 

There is another approach to find the expression for 
correlation functions. It was shown in \cite{JMN} that 
the correlation functions of the eight-vertex model satisfy 
a set of difference equations, quantum Knizhnik--Zamolodchikov 
equation of level $-4$. On the basis of difference equation approach, 
we obtained the expression of the spontaneous polarization of 
$(\mathbb{Z}/n\mathbb{Z})$-symmetric model \cite{SPn}. 
In this paper we show that the expressions for 
the spontaneous polarization of 
$(\mathbb{Z}/n\mathbb{Z})$-symmetric model obtained 
on the basis of free field representation formalism reproduces 
the known result in \cite{SPn}. 
This coincidence indicates 
the relevance of the free field representation formalism. 

The present paper is organized as follows. In section 2 
we review the basic definitions of 
$(\mathbb{Z}/n\mathbb{Z})$-symmetric model \cite{Bela}, 
the corresponding dual face model \cite{JMO}, 
and the vertex-face correspondence. In section 3 we introduce 
the CTM(corner transfer matrix) Hamiltonians and the vertex operators 
of $(\mathbb{Z}/n\mathbb{Z})$-symmetric model and 
$A^{(1)}_{n-1}$ model, and also introduce the tail operators 
which relates those two CTM Hamiltonians. In section 4 we construct 
the free field formalism of $(\mathbb{Z}/n\mathbb{Z})$-symmetric model. 
In section 5 we present trace formulae for correlation functions of 
$(\mathbb{Z}/n\mathbb{Z})$-symmetric model. Furthermore, 
we calculate the spontaneous polarization of 
$(\mathbb{Z}/n\mathbb{Z})$-symmetric model in this formalism. 
Sections 4 and 5 are main original parts of the present paper. 
In section 6 we give some concluding remarks. 

\section{Basic definitions}

The present section aims to formulate 
the problem, thereby fixing the notation. 

\subsection{Theta functions} 

Jacobi theta function with two pseudo-periods $1$ and 
$\tau$\,(${\rm Im}\,\tau >0$) are defined as follows: 
\begin{equation}
\vartheta\left[\begin{array}{c} a \\ b \end{array} \right]
(v;\tau ): =\displaystyle\sum_{m\in \mathbb{Z}} 
\exp \left\{ \pi \sqrt{-1}(m+a)~\left[ (m+a)\tau 
+2(v+b) \right] \right\}, \label{Rieth}
\end{equation}
for $a,b\in\mathbb{R}$. 
Let $n\in\mathbb{Z}_{\geqslant 2}$ and 
$r\in\mathbb{R}$ such that $r> n-1$, and also fix 
the parameter $x$ such that $0<x<1$. 
We will use the abbreviations, 
\begin{equation}
[v]=x^{\frac{v^2}{r}-v}\Theta_{x^{2r}}(x^{2v}), 
~~~~
[v]'=x^{\frac{v^2}{r-1}-v}\Theta_{x^{2r-2}}(x^{2v}), 
\end{equation}
where 
\begin{eqnarray*}
&&\Theta_{q}(z)=(z; q)_\infty 
(qz^{-1}; q)_\infty (q; q)_\infty =
\sum_{m\in\mathbb{Z}} q^{m(m-1)/2}(-z)^m, \\
&&(z; q_1 , \cdots , q_m )= 
\prod_{i_1 , \cdots , i_m \geqslant 0} 
(1-zq_1^{i_1} \cdots q_m^{i_m}). 
\end{eqnarray*}
Note that 
$$
\vartheta\left[\begin{array}{c} 1/2 \\ -1/2 \end{array} \right]
\left( \frac{v}{r}, \frac{\pi\sqrt{-1}}{\epsilon r} \right)
=\sqrt{\frac{\epsilon r}{\pi}}
\exp\,\left(-\frac{\epsilon r}{4}\right)[v], 
$$
where $x=e^{-\epsilon}$ ($\epsilon >0$). 

For later conveniences we also 
introduce the following symbols
\begin{eqnarray}
r_{l}(v)&=&z^{\frac{r-1}{r}\frac{n-l}{n}}
\frac{g_{l}(z^{-1})}{g_{l}(z)}, ~~~~
g_{l}(z)=
\frac{\{x^{2n+2r-l-1}z\}
\{x^{l+1}z\}}
{\{x^{2n-l+1}z\}\{x^{2r+l-1}z\}}, 
\end{eqnarray}
where $z=x^{2v}$, $1\leqslant l\leqslant n$ and 
\begin{equation}
\{z\}=(z;x^{2r},x^{2n})_\infty . 
\label{eq:{z}}
\end{equation}
These factors will appear in the commutation relations 
among the type I vertex operators. 

The integral kernel for the type I 
vertex operators will be given as the products of 
the following elliptic functions 
\begin{eqnarray}
f(v,w)=\frac{[v+\frac{1}{2}-w]}{[v-\frac{1}{2}]}, &~~~&
g(v)=\frac{[v-1]}{[v+1]}. 
\end{eqnarray}

\subsection{Belavin's vertex model}

Let $V=\mathbb{C}^n$ and 
$\{ \varepsilon _\mu \}_{0 \leqslant \mu \leqslant n-1}$ be 
the standard orthonormal basis with the inner 
product $\langle \varepsilon _\mu , 
\varepsilon _\nu \rangle =\delta_{\mu \nu}$. 
Belavin's $(\mathbb{Z}/n\mathbb{Z})$-symmetric model is 
a vertex model on a two-dimensional square lattice ${\cal L}$ 
such that the state variables take on values of 
$(\mathbb{Z}/n\mathbb{Z})$-spin. 
In the original papers \cite{Bela,RT}, 
the $R$-matrix in the disordered phase is given. 
For the present purpose, we need the following $R$-matrix: 
\begin{equation}
R(v)=\dfrac{[1]}{[1-v]}r_1 (v)\overline{R}(v), ~~~~
\overline{R}(v)=\dfrac{1}{n}
\sum_{\mbox{\footnotesize\boldmath $\alpha $}\in G_n} 
\dfrac{\vartheta \left[\begin{array}{c} \frac{1}{2}-\frac{\alpha_1}{n} \\ 
\frac{1}{2}+\frac{\alpha_2}{n} \end{array} \right]
\left( \dfrac{1}{nr}-\dfrac{v}{r} ; 
\dfrac{\pi\sqrt{-1}}{\epsilon r} \right)}{
\vartheta \left[\begin{array}{c} \frac{1}{2}-\frac{\alpha_1}{n} \\ 
\frac{1}{2}+\frac{\alpha_2}{n} \end{array} \right]\left( \dfrac{1}{nr}; 
\dfrac{\pi\sqrt{-1}}{\epsilon r} \right)}
I_{\mbox{\footnotesize\boldmath $\alpha $}} 
\otimes I_{\mbox{\footnotesize\boldmath $\alpha $}}^{-1}. 
\label{eq:Bel-sol}
\end{equation}
Here $G_n =(\mathbb{Z}/n\mathbb{Z}) \times 
(\mathbb{Z}/n\mathbb{Z})$, and 
$I_{\mbox{\footnotesize\boldmath $\alpha $}}=
g^{\alpha_1}h^{\alpha_2}$ for 
$\mbox{\boldmath $\alpha $}=
(\alpha_1 , \alpha_2 )$, where 
\begin{equation}
gv_i =\omega ^i v_i ,~~~~~~
hv_i =v_{i-1}, \label{gh}
\end{equation}
with $\omega =\exp (2\pi \sqrt{-1}/n)$. 
We assume that the parameters $v$, $\epsilon$ and $r$ lie 
in the so-called principal regime: 
\begin{equation}
\epsilon >0, ~~ r>1, ~~ 0<v<1. 
\label{eq:principal}
\end{equation}
When $n=2$ the principal regime (\ref{eq:principal}) 
lies in one of the antiferroelectric phases of 
the eight-vertex model \cite{ESM}. We describe 
$n$ kinds of ground states of $(\mathbb{Z}/n\mathbb{Z})$-symmetric 
model in the principal regime in section 3.1. 

The $R$-matrix satisfies the Yang-Baxter equation (YBE) 
\begin{eqnarray}
R_{12}(v_1 -v_2 )
R_{13}(v_1 -v_3 )
R_{23}(v_2 -v_3 )=
R_{23}(v_2 -v_3 )
R_{13}(v_1 -v_3 )
R_{12}(v_1 -v_2 ), \label{YBE}
\end{eqnarray}
where $R_{ij}(v)$ denotes the matrix on $V^{\otimes 3}$, 
which acts as $R(v)$ on the $i$-th and $j$-th components and 
as identity on the other one. 

If $i+k=j+l$ (mod $n$), 
the elements of $R$-matirix $\overline{R}(v)^{ik}_{jl}$ 
is given as follows: 
\begin{equation}
\overline{R}(v)^{ik}_{jl}=
\displaystyle\frac{h(v)
\vartheta \left[\begin{array}{c} \frac{1}{2} \\ 
\frac{1}{2}+\frac{k-i}{n} \end{array} \right]
\left( \dfrac{1-v}{nr} ; 
\dfrac{\pi\sqrt{-1}}{n\epsilon r} \right)}
{\vartheta \left[\begin{array}{c} \frac{1}{2} \\ 
\frac{1}{2}+\frac{j-k}{n} \end{array} \right]
\left( \dfrac{v}{nr} ; 
\dfrac{\pi\sqrt{-1}}{n\epsilon r} \right)
\vartheta \left[\begin{array}{c} \frac{1}{2} \\ 
\frac{1}{2}+\frac{j-i}{n} \end{array} \right]
\left( \dfrac{1}{nr} ; 
\dfrac{\pi\sqrt{-1}}{n\epsilon r} \right)}, 
\label{Bel}
\end{equation}
where 
$$
h(v)=\prod_{j=0}^{n-1} 
\vartheta \left[\begin{array}{c} \frac{1}{2} \\ 
\frac{1}{2}+\frac{j}{n} \end{array} \right]
\left( \frac{v}{nr} ; 
\frac{\pi\sqrt{-1}}{n\epsilon r} \right)\left/
\;\prod_{j=1}^{n-1} \vartheta \left[\begin{array}{c} \frac{1}{2} \\ 
\frac{1}{2}+\frac{j}{n} \end{array} \right]\right.
\left( 0 ; 
\frac{\pi\sqrt{-1}}{n\epsilon r} \right), 
$$
and otherwize $\overline{R}(v)^{ik}_{jl}=0$. 
\label{prop:Bel}

Note that the weights (\ref{Bel}) reproduce those of the eight-vertex 
model in the principal regime when $n=2$ \cite{ESM}. 

\subsection{The weight lattice and the root lattice 
of $A^{(1)}_{n-1}$}

Let $V=\mathbb{C}^n$ and 
$\{ \varepsilon _\mu \}_{0 \leqslant \mu \leqslant n-1}$ be 
the standard orthonormal basis as before. 
The weight lattice of $A^{(1)}_{n-1}$ 
is defined as follows: 
\begin{equation}
P=\bigoplus_{\mu =0}^{n-1} 
\mathbb{Z} \bar{\varepsilon}_\mu , 
\label{eq:wt-lattice}
\end{equation}
where 
$$
\bar{\varepsilon}_\mu =
\varepsilon _\mu -\varepsilon , 
~~~~\varepsilon =\frac{1}{n}\sum_{\mu =0}^{n-1} 
\varepsilon _\mu . 
$$
We denote the fundamental weights 
by $\omega_\mu\,(1\leqslant \mu \leqslant n-1)$ 
$$
\omega_\mu =\sum_{\nu =0}^{\mu -1}
\bar{\varepsilon }_\nu , 
$$
and also denote the simple roots by 
$\alpha_\mu \,(1\leqslant \mu \leqslant n-1)$ 
$$
\alpha_\mu =\varepsilon_{\mu-1} -
\varepsilon _{\mu}=\bar{\varepsilon}_{\mu -1}
-\bar{\varepsilon}_{\mu}. 
$$
The root lattice of $A^{(1)}_{n-1}$ 
is defined as follows: 
\begin{equation}
Q=\bigoplus_{\mu =1}^{n-1} 
\mathbb{Z} \alpha_\mu , 
\label{eq:rt-lattice}
\end{equation}

For $a\in P$ we set 
\begin{equation}
a_{\mu\nu}=\bar{a}_\mu-\bar{a}_\nu , ~~~~ 
\bar{a}_\mu =\langle a+\rho , 
\varepsilon_\mu \rangle =\langle a+\rho , 
\bar{\varepsilon}_\mu \rangle , ~~~~ 
\rho =\sum_{\mu =1}^{n-1} \omega_\mu . 
\end{equation}

Useful formulae are: 
$$
\langle \bar{\varepsilon}_{\mu}, \varepsilon_{\nu} \rangle 
=\langle \bar{\varepsilon}_{\mu}, \bar{\varepsilon}_{\nu} \rangle 
=\delta_{\mu\nu}-\dfrac{1}{n}, ~~~~
\langle \alpha_\mu , \omega_\nu \rangle =\delta_{\mu\nu}, 
$$
$$
\langle \bar{\varepsilon}_{\mu}, \omega_{\nu} \rangle 
=\theta (\mu <\nu)-\dfrac{\nu}{n}, ~~~~
\langle \omega_{\mu}, \omega_{\nu} \rangle 
=\mbox{min}\,(\mu , \nu )-\dfrac{\mu\nu}{n}. 
$$
When $a+\rho =\displaystyle\sum_{\mu=0}^{n-1} k^\mu \omega_\mu$, 
we have $a_{\mu\nu}=k^{\mu +1}+\cdots +k^\nu$ when 
$\mu <\nu$, and 
$$
\langle a+\rho , a+\rho \rangle =
\dfrac{1}{n}\sum_{\mu <\nu} a_{\mu\nu}^2, ~~~~
\langle a+\rho , \rho \rangle =
\dfrac{1}{2}\sum_{\mu <\nu} a_{\mu\nu}. 
$$
Let $\displaystyle\sum_{\mu=0}^{n-1} k^\mu =r$, where 
$a+\rho =\displaystyle\sum_{\mu=0}^{n-1} k^\mu \omega_\mu$, 
then we denote $a\in P_{r-n}$. 

\subsection{The $A^{(1)}_{n-1}$ face model}

An ordered pair $(a,b) \in P^2_{r-n}$ 
is called {\it admissible} if $b=a+\bar{\varepsilon}_\mu$, 
for a certain $\mu\,(0\leqslant \mu \leqslant n-1)$. 
For $(a, b, c, d)\in P^4_{r-n}$ let 
$\displaystyle W 
\left[ \left. \begin{array}{cc} 
c & d \\ b & a \end{array} 
\right| v \right] $ 
be the Boltzmann weight of the 
$A^{(1)}_{n-1}$ model for the state configuration 
$\displaystyle 
\left[ \begin{array}{cc} 
c & d \\ b & a \end{array} \right] $ 
round a face. 
Here the four states $a, b, c$ and $d$ are 
ordered clockwise from the SE corner. 
In this model $W 
\left[ \left. \begin{array}{cc} 
c & d \\ b & a \end{array} \right| 
v \right] =0~~$ 
unless the four pairs $(a,b), (a,d), (b,c)$ 
and $(d,c)$ are admissible. 
Non-zero Boltzmann weights are parametrized in terms of 
the elliptic theta function of the spectral parameter $v$ 
as follows: 
\begin{equation}
\begin{array}{rcl}
W
\left[ \left. \begin{array}{cc} 
a + 2 \bar{\varepsilon }_\mu & a+\bar{\varepsilon }_\mu \\ 
a+\bar{\varepsilon }_\mu & a \end{array} \right| 
v  \right] 
& = & r_1 (v), \\
~ & ~ & ~ \\
W
\left[ \left. \begin{array}{cc} 
a+\bar{\varepsilon }_\mu +\bar{\varepsilon }_\nu 
& a+\bar{\varepsilon }_\mu \\ 
a+\bar{\varepsilon }_\nu & a 
\end{array} \right| v \right] & = & -r_1 (v)
\dfrac{[v][a_{\mu\nu}+1]}{[1-v][a_{\mu\nu}]} 
~~~~(\mu \neq \nu ), \label{BW} \\
~ & ~ & ~ \\
W
\left[ \left. \begin{array}{cc} 
a+\bar{\varepsilon }_\mu +\bar{\varepsilon }_\nu 
& a+\bar{\varepsilon }_\mu \\ 
a+\bar{\varepsilon }_\mu & a 
\end{array} \right| v \right] & = & r_1 (v)
\dfrac{[1][v+a_{\mu\nu}]}{[1-v][a_{\mu\nu}]} 
~~~~(\mu\neq \nu). 
\end{array}
\end{equation}
We consider so-called Regime III in the model, i.e., 
$0<v<1$. 

The Boltzmann weights (\ref{BW}) 
solve the Yang-Baxter equation 
for the face model \cite{JMO}: 
\begin{equation}
\begin{array}{cc}
~ & 
\displaystyle \sum_{g} 
W\left[ \left. 
\begin{array}{cc} d & e \\ c & g \end{array} \right| 
v_1 \right]
W\left[ \left. 
\begin{array}{cc} c & g \\ b & a \end{array} \right| 
v_2 \right]
W\left[ \left. 
\begin{array}{cc} e & f \\ g & a \end{array} \right| 
v_1 -v_2 \right] \\
~ & ~ \\
= & \displaystyle \sum_{g} 
W\left[ \left. 
\begin{array}{cc} g & f \\ b & a \end{array} \right| 
v_1 \right]
W\left[ \left. 
\begin{array}{cc} d & e \\ g & f \end{array} \right| 
v_2 \right]
W\left[ \left. 
\begin{array}{cc} d & g \\ c & b \end{array} \right| 
v_1 -v_2 \right]
\end{array} \label{STR}
\end{equation}

\subsection{Vertex-face correspondence}

In this paper we use the $R$-matrix of 
$(\mathbb{Z}/n\mathbb{Z})$-symmetric model in 
the principal regime while Belavin's 
original paper used the one in the disordered phase. 
Thus, we need different intertwining vectors 
from the one by Jimbo-Miwa-Okado \cite{JMO}. 

Let 
\begin{equation}
t(v)^a_{a-\bar{\varepsilon}_\mu}=
\sum_{\nu =0}^{n-1} \varepsilon_\nu 
\vartheta  \left[\begin{array}{c} 0 \\ 
\frac{1}{2}+\frac{\nu}{n} \end{array} \right] 
\left( \frac{v}{nr}+\frac{\bar{a}_\mu}{r}; 
\frac{\pi\sqrt{-1}}{n\epsilon r} \right). 
\label{eq:int-vec}
\end{equation}
Then we have (cf. figure 1) 
\begin{equation}
R(v_1-v_2)t (v_1)_a^d\otimes t (v_2)_d^c=
\sum_{b} t(v_1)_b^c \otimes t (v_2)_a^b 
W\left[ \left. 
\begin{array}{cc} c & d \\ b & a \end{array} \right| 
v_1 -v_2 \right]. 
\label{eq:Rtt=Wtt}
\end{equation}

\unitlength 1mm
\begin{picture}(100,20)
\put(23,0){
\begin{picture}(101,0)
\put(20,3){\begin{picture}(101,0)
\put(10,10){\vector(-1,0){10}}
\put(10,0){\vector(0,1){10}}
\put(8.8,5){\vector(-1,0){10}}
\put(8.,4.4){\scriptsize{$<$}}
\put(-4.5,4.2){$v_1$}
\put(4.15,8.5){\scriptsize{$\vee$}}
\put(5,8.8){\vector(0,-1){10}}
\put(3.9,-3.8){$v_2$}
\put(-2.5,10.5){$c$}
\put(10.5,-1.5){$a$}
\put(10.5,10.1){$d$}
\put(17,4){$=\;\displaystyle\sum_{b}$} 
\end{picture}
}
\put(58,3){\begin{picture}(101,0)
\put(10,0){\vector(-1,0){10}}
\put(10,0){\vector(0,1){10}}
\put(0,0){\vector(0,1){10}}
\put(10,10){\vector(-1,0){10}}
\put(-2.,4.4){\scriptsize{$<$}}
\multiput(0,5)(2.2,0){6}{\line(1,0){1.2}}
\put(-1.2,5){\vector(-1,0){2.5}}
\put(-7,4.2){$v_1$}
\put(4.15,-1.5){\scriptsize{$\vee$}}
\multiput(5,0)(0,2.2){6}{\line(0,1){1.2}}
\put(5,-1.2){\vector(0,-1){2.5}}
\put(3.9,-5.8){$v_2$}
\put(10.5,10.1){$d$}
\put(-2.5,10.5){$c$}
\put(10.5,-1.5){$a$}
\put(-2.5,-1.8){$b$}
\end{picture}
}
\end{picture}
}
\end{picture}

\vspace{2mm}

\begin{center}
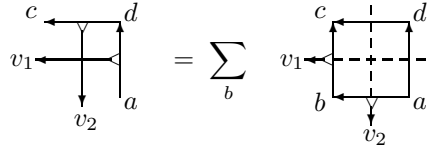
Figure 1. Picture representation of vertex-face 
correspondence. 
\end{center}

\section{Vertex-face transformation}

The basic objects in the vertex operator approach are 
the corner transfer matrices (CTM) and the vertex operators \cite{JMbk}. 
In subsections 3.1 and 3.2 we recall the CTM Hamiltonians, the type I 
vertex operators and the space of states of 
$(\mathbb{Z}/n\mathbb{Z})$-symmetric model and the $A^{(1)}_{n-1}$ 
model, respectively. 

In \cite{LaP} Lashkevich and Pugai introduced the nonlocal operator 
called the tail operator, in order to express correlation functions of the 
eight-vertex model in terms of those of the SOS model. 
In subsection 3.3 we introduce the tail operator for the present purpose; 
i.e., in order to express correlation functions of the 
$(\mathbb{Z}/n\mathbb{Z})$-symmetric model in terms of those of the 
$A^{(1)}_{n-1}$ model. The commutation relations among the tail operators 
and the type I vertex operators are given in subsection 3.4. 

\subsection{CTM Hamiltonian for the vertex model}

Let us consider the `low temperature' limit $x\rightarrow 0$. 
Then the elements of $R$-matrix behave as 
\begin{equation}
R^{\mu\nu}_{\mu'\nu'}(v)\sim \zeta^{H_v(\mu ,\nu )}\delta^\mu_{\nu'} 
\delta^\nu_{\mu'}, 
\label{eq:B-lowtemp}
\end{equation}
where $z=x^{2v}=\zeta^n$ and 
\begin{equation}
H_v (\mu ,\nu )=\left\{ \begin{array}{ll} 
\mu -\nu -1 & \mbox{if $0\leqslant \nu <\mu\leqslant n-1$} \\
n-1+\mu -\nu & \mbox{if $0\leqslant \mu\leqslant \nu\leqslant n-1$}
\end{array} \right. 
\label{eq:H_v}
\end{equation}
Thus the CTM(corner transfer matrix) 
Hamiltonian of $(\mathbb{Z}/n\mathbb{Z})$-symmetric 
model in the principal regime is given as follows: 
\begin{equation}
H_{CTM}(\mu_1 , \mu_2 , \mu_3, \cdots )=
\sum_{j=1}^\infty jH_v (\mu_j , \mu_{j+1}). 
\end{equation}
The CTM Hamiltonian diverges unless $\mu_j =i+1-j$ (mod $n$) for 
$j\gg 0$ and a certain $0\leqslant i\leqslant n-1$. 

Let ${\cal H}^{(i)}$ be the $\mathbb{C}$-vector space 
spanned by the half-infinite pure tensor vectors of the forms\footnote{ 
We fix ${\cal H}^{(i)}$ by (\ref{eq:H^i}) such that it coincides with 
$V(\omega_{i})$, the level $1$ highest weight irreducible 
$U_q (\widehat{\mathfrak{s}\mathfrak{l}_n})$-module, 
in the trigonometric limit 
$r\rightarrow \infty$. For example, see \cite{Koy}, keeping in mind 
that our $i$ should be read as $-i$ in \cite{Koy}. }
\begin{equation}
\varepsilon_{\mu_1}\otimes \varepsilon_{\mu_2}\otimes 
\varepsilon_{\mu_2}\otimes \cdots 
~~~~ \mbox{with $\mu_j\in \mathbb{Z}/n\mathbb{Z}$, 
$\mu_j=i+1-j$ (mod $n$) for $j\gg 0$}. 
\label{eq:H^i}
\end{equation}
Let ${\cal H}^{*(i)}$ be the dual of ${\cal H}^{(i)}$ 
spanned by the half-infinite pure tensor vectors of the forms 
\begin{equation}
\cdots \otimes \varepsilon_{\mu_{-2}}\otimes \varepsilon_{\mu_{-1}}\otimes 
\varepsilon_{\mu_0}
~~~~ \mbox{with $\mu_j\in \mathbb{Z}/n\mathbb{Z}$, 
$\mu_j=i+1-j$ (mod $n$) for $j\ll 0$}. 
\end{equation}

Introduce the type I vertex operator by the following 
half-infinite transfer matrix 
\begin{equation}
\Phi^\mu (v_1 -v_2)=
\unitlength 0.5mm
\begin{picture}(100,20)
\put(15,-18){\begin{picture}(100,0)
\put(60,20){\vector(-1,0){60}}
\put(10,30){\vector(0,-1){20}}
\put(20,30){\vector(0,-1){20}}
\put(30,30){\vector(0,-1){20}}
\put(40,30){\vector(0,-1){20}}
\put(4,22){$\mu$}
\put(-7,18){$v_1$}
\put(9,5){$v_2$}
\put(19,5){$v_2$}
\put(29,5){$v_2$}
\put(39,5){$v_2$}
\put(45,22){$\cdots$}
\end{picture}
}
\end{picture}
\label{eq:B-I}
\end{equation}

~

\noindent Then the operator (\ref{eq:B-I}) is an intertwiner 
from ${\cal H}^{(i)}$ to ${\cal H}^{(i+1)}$. 
The type I vertex operators satisfy the following 
commutation relation: 
\begin{equation}
\Phi^\mu (v_1)\Phi^\nu (v_2)=
\sum_{\mu',\nu'} R(v_1-v_2)^{\mu\nu}_{\mu'\nu'} 
\Phi^{\nu'} (v_2)\Phi^{\mu'} (v_1). 
\label{eq:RPhiPhi}
\end{equation}

Introduce the corner transfer matrix in the SE(South-East) corner 

\unitlength 1mm
\begin{picture}(140,65)(-5,-5)
\put(0,30){$A^{(i)}_{SE}(v_{1}-v_{2})^{\mu_1\mu_2\mu_3\mu_4\cdots}
_{\mu'_1\mu'_2\mu'_3\mu'_4\cdots}=$}
\put(46,0){\begin{picture}(100,0)
\put(20,10){\vector(-1,0){20}}
\put(30,20){\vector(-1,0){30}}
\put(40,30){\vector(-1,0){40}}
\put(50,40){\vector(-1,0){50}}
\put(10,50){\vector(0,-1){50}}
\put(20,50){\vector(0,-1){40}}
\put(30,50){\vector(0,-1){30}}
\put(40,50){\vector(0,-1){20}}
\put(-5,9){$v_{2}$}
\put(-5,19){$v_{2}$}
\put(-5,29){$v_{2}$}
\put(-5,39){$v_{2}$}
\put(39,52){$v_{1}$}
\put(29,52){$v_{1}$}
\put(19,52){$v_{1}$}
\put(9,52){$v_{1}$}
\put(4,41){$\mu_1$}
\put(4,31){$\mu_2$}
\put(4,21){$\mu_3$}
\put(4,11){$\mu_4$}
\put(4.5,3){$\vdots$}
\put(11,44){$\mu'_1$}
\put(21,44){$\mu'_2$}
\put(31,44){$\mu'_3$}
\put(41,44){$\mu'_4$}
\put(48,44){$\cdots$}
\put(31,15.3){$\ddots$}
\end{picture}
}
\end{picture}

\noindent The diagonal form of $A^{(i)}_{SE}(v)$ can be determined from the 
`low temperature' limit of the $R$-matrix (\ref{eq:B-lowtemp}--\ref{eq:H_v}): 
\begin{equation}
A^{(i)}_{SE}(v) \sim \zeta^{H_{CTM}}=z^{\frac{1}{n}H_{CTM}}: 
{\cal H}^{(i)}\rightarrow {\cal H}^{(i)}, 
\end{equation}
where $\sim$ refers to an equality modulo a divergent scalar in the infinite 
lattice limit. 
Likewise other three types of the corner transfer matrices 
are given as follows: 
\begin{equation}
\begin{array}{ll}
A^{(i)}_{NE}(v): & {\cal H}^{(i)}\rightarrow {\cal H}^{*(i)}, \\
A^{(i)}_{NW}(v): & {\cal H}^{*(i)}\rightarrow {\cal H}^{*(i)}, \\
A^{(i)}_{SW}(v): & {\cal H}^{*(i)}\rightarrow {\cal H}^{(i)}, \\
\end{array}
\end{equation}
where NE, NW and SW stand for the corners North-East, North-West 
and South-West. 
It seems to be rather general \cite{ESM} that the product of four CTMs 
in the infinite lattice limit is independent of $v$: 
\begin{equation}
\rho^{(i)}=A^{(i)}_{SE}(v)A^{(i)}_{SW}(v)
A^{(i)}_{NW}(v)A^{(i)}_{NE}(v)
=x^{2H_{CTM}}. 
\end{equation}
Since $H(\mu_j , \mu_{j+1})$ takes on value of 
$\{ 0,1,\cdots , n-1\}$, the eigenvalues of 
$H_{CTM}$ are of the form 
$$
N=\sum_{j=1}^\infty jm_j , ~~~~ 0\leqslant m_j \leqslant n-1. 
$$
This stands for the partition of $N$ 
such that the multiplicity of each $j$ is 
at most $n-1$. Thus, 
the character is given by 
\begin{equation}
\chi^{(i)}=\mbox{tr}_{{\cal H}^{(i)}}\,(\rho^{(i)})
=\dfrac{(x^{2n};x^{2n})_\infty}{(x^2;x^2)_\infty}. 
\label{eq:chi^i}
\end{equation}

\subsection{CTM for the $A^{(1)}_{n-1}$ model}

After gauge transformation \cite{JMO}, 
the CTM Hamiltonian of $A^{(1)}_{n-1}$ model 
in regime III is given as follows: 
\begin{equation}
\begin{array}{rcl}
H_{CTM}(a_0 , a_1 , a_2, \cdots )&=&
\displaystyle\sum_{j=1}^\infty jH_f (a_{j-1} , a_{j}, a_{j+1}), \\
H_f(a+\bar{\varepsilon }_\mu +\bar{\varepsilon }_\nu , 
a+\bar{\varepsilon }_\mu , a )&=&
\dfrac{1}{n}H_v (\nu , \mu ), \end{array}
\label{eq:lowtemp-face}
\end{equation}
where $H_v (\nu , \mu )$ is defined by (\ref{eq:H_v}). 
The CTM Hamiltonian diverges unless $a_j =\xi +\omega_{i+1-j}$ 
for $j\gg 0$ and a certain $\xi\in P_{r-n-1}$ 
and $0\leqslant i\leqslant n-1$. 

For $k=a+\rho , l=\xi +\rho$ and $0\leqslant i\leqslant n-1$, 
let ${\cal H}^{(i)}_{l,k}$ be the space of admissible paths 
$(a_0 , a_1, a_2, \cdots )$ such that 
\begin{equation}
a_0 =a, ~~~ a_{j} -a_{j+1}\in \left\{ 
\bar{\varepsilon}_0 , \bar{\varepsilon}_1 , 
\cdots , \bar{\varepsilon}_{n-1} 
\right\}, \mbox{ for $j=1, 2, 3, \cdots$, }~~~~ 
a_j=\xi +\omega_{i+1-j} \mbox{ for 
$j\gg 0$}. 
\end{equation}
Also, let ${\cal H}^{*(i)}_{l,k}$ be the space of admissible paths 
$(\cdots , a_{-2} , a_{-1}, a_{0})$ such that 
\begin{equation}
a_0 =a, ~~~ a_{j} -a_{j+1}\in \left\{ 
\bar{\varepsilon}_0 , \bar{\varepsilon}_1 , 
\cdots , \bar{\varepsilon}_{n-1} 
\right\}, \mbox{ for $j=1, 2, 3, \cdots$, }~~~~ 
a_j=\xi +\omega_{i+1-j} \mbox{ for 
$j\ll 0$}. 
\end{equation}
Introduce the type I vertex operator by the following 
half-infinite transfer matrix 
\begin{equation}
\hspace{2cm} \Phi(v_1 -v_2)_a^{a+\bar{\varepsilon}_\mu}=
\unitlength 1mm
\begin{picture}(100,10)
\put(-10,-60){\begin{picture}(101,0)
\put(20,55){\vector(0,1){10}}
\put(30,55){\vector(0,1){10}}
\put(30,55){\vector(-1,0){10}}
\put(30,65){\vector(-1,0){10}}
\put(40,55){\vector(0,1){10}}
\put(40,55){\vector(-1,0){10}}
\put(40,65){\vector(-1,0){10}}
\put(50,55){\vector(0,1){10}}
\put(50,55){\vector(-1,0){10}}
\put(50,65){\vector(-1,0){10}}
\put(60,55){\vector(-1,0){10}}
\put(60,65){\vector(-1,0){10}}
\put(19,53){$a$}
\put(16,67){$a\!+\!\bar{\varepsilon}_\mu$}
\multiput(18,60)(2,0){23}{\line(1,0){1}}
\put(17,60){\vector(-1,0){1.5}}
\multiput(25,53)(0,2){8}{\line(0,1){1}}
\multiput(35,53)(0,2){8}{\line(0,1){1}}
\multiput(45,53)(0,2){8}{\line(0,1){1}}
\put(25,52){\vector(0,-1){1.5}}
\put(35,52){\vector(0,-1){1.5}}
\put(45,52){\vector(0,-1){1.5}}
\put(24,48){$v_2$}
\put(34,48){$v_2$}
\put(44,48){$v_2$}
\put(12,59.5){$v_1$}
\end{picture}
}
\end{picture}
\label{eq:F-I}
\end{equation}

\vspace{7mm}

\noindent Then the operator (\ref{eq:F-I}) is an intertwiner 
from ${\cal H}^{(i)}_{l,k}$ to 
${\cal H}^{(i+1)}_{l,k+\bar{\varepsilon}_\mu}$. 
The type I vertex operators satisfy the following 
commutation relation: 
\begin{equation}
\Phi (v_1)^c_b\Phi (v_2)^b_a=
\sum_{d} W\left[ \left. \begin{array}{cc} 
c & d \\ 
b & a \end{array} \right| v_1-v_2 \right]
\Phi (v_2)^{c}_d\Phi (v_1)^{d}_a . 
\label{eq:Wphiphi}
\end{equation}

Introduce the corner transfer matrix of the $A^{(1)}_{n-1}$ model 
in the SE corner 

\unitlength 1mm
\begin{picture}(140,65)(-15,-5)
\put(-5,30){$A^{(l,k)}_{SE}(v_{1}-v_{2})^{a_0a_1a_2a_3a_4\cdots}
_{a_0a'_1a'_2a'_3a'_4\cdots}=$}
\put(46,0){\begin{picture}(100,0)
\multiput(10,2)(0,2){24}{\line(0,1){1}}
\multiput(20,12)(0,2){19}{\line(0,1){1}}
\multiput(30,22)(0,2){14}{\line(0,1){1}}
\multiput(40,32)(0,2){9}{\line(0,1){1}}
\multiput(2,40)(2,0){24}{\line(1,0){1}}
\multiput(2,30)(2,0){19}{\line(1,0){1}}
\multiput(2,20)(2,0){14}{\line(1,0){1}}
\multiput(2,10)(2,0){9}{\line(1,0){1}}
\put(0,40){\vector(-1,0){1.5}}
\put(0,30){\vector(-1,0){1.5}}
\put(0,20){\vector(-1,0){1.5}}
\put(0,10){\vector(-1,0){1.5}}
\put(10,1.5){\vector(0,-1){1.5}}
\put(20,11.5){\vector(0,-1){1.5}}
\put(30,21.5){\vector(0,-1){1.5}}
\put(40,31.5){\vector(0,-1){1.5}}
\put(-5,9){$v_{1}$}
\put(-5,19){$v_{1}$}
\put(-5,29){$v_{1}$}
\put(-5,39){$v_{1}$}
\put(39,50){$v_{2}$}
\put(29,50){$v_{2}$}
\put(19,50){$v_{2}$}
\put(9,50){$v_{2}$}
\put(5,0){\line(0,1){45}}
\put(15,5){\line(0,1){40}}
\put(25,15){\line(0,1){30}}
\put(35,25){\line(0,1){20}}
\put(45,35){\line(0,1){10}}
\put(5,5){\line(1,0){10}}
\put(5,15){\line(1,0){20}}
\put(5,25){\line(1,0){30}}
\put(5,35){\line(1,0){40}}
\put(5,45){\line(1,0){45}}
\put(1,46){$a_0$}
\put(1,34){$a_1$}
\put(1,24){$a_2$}
\put(1,14){$a_3$}
\put(1,4){$a_4$}
\put(1.5,-2){$\vdots$}
\put(13,47){$a'_1$}
\put(23,47){$a'_2$}
\put(33,47){$a'_3$}
\put(43,47){$a'_4$}
\put(50,47){$\cdots$}
\put(31,15.3){$\ddots$}
\end{picture}
}
\end{picture}

\noindent The diagonal form of $A^{(l,k)}_{SE}(v)$ can be determined from the 
`low temperature' limit (\ref{eq:lowtemp-face}): 
\begin{equation}
A^{(l,k)}_{SE}(v) \sim \zeta^{H_{CTM}}=z^{\frac{1}{n}H_{CTM}}: 
{\cal H}^{(i)}_{l,k}\rightarrow {\cal H}^{(i)}_{l,k}, 
\end{equation}
where $\sim$ refers to an equality modulo a divergent scalar in the infinite 
lattice limit. 
Likewise other three types of the corner transfer matrices 
are given as follows: 
\begin{equation}
\begin{array}{ll}
A^{(l,k)}_{NE}(v): & {\cal H}^{(i)}_{l,k}\rightarrow {\cal H}^{*(i)}_{l,k}, \\
A^{(l,k)}_{NW}(v): & {\cal H}^{*(i)}_{l,k}\rightarrow {\cal H}^{*(i)}_{l,k}, \\
A^{(l,k)}_{SW}(v): & {\cal H}^{*(i)}_{l,k}\rightarrow {\cal H}^{(i)}_{l,k}. \\
\end{array}
\end{equation}
The product of four CTMs for $A^{(1)}_{n-1}$ model 
in the infinite lattice limit is also independent of $v$ \cite{JMO}: 
\begin{equation}
\rho^{(i)}_{l,k}=G_a x^{2nH_{l,k}^{(i)}}, 
\label{eq:rho_lk}
\end{equation}
where 
$$
G_a =\prod_{\mu <\nu} [a_{\mu\nu}]. 
$$
The character of $A^{(1)}_{n-1}$ model 
was obtained in \cite{JMO}: 
\begin{equation}
\chi^{(i)}_{l,k}=\mbox{tr}_{{\cal H}^{(i)}_{l,k}}\,(\rho^{(i)}_{l,k})
=\dfrac{x^{n |\beta_1 k+\beta_2l|^2}}
{(x^{2n};x^{2n})^{n-1}_\infty}G_a, 
\label{eq:chi_lk}
\end{equation}
where 
\begin{equation}
t^2 -\beta_0 t-1=(t-\beta_1)(t-\beta_2), ~~~~ 
\beta_0 =\dfrac{1}{\sqrt{r(r-1)}}, ~~~~ \beta_1<\beta_2. 
\label{eq:beta_12}
\end{equation}

We notice the following sum formula: 
\begin{equation}
\sum_{k\equiv l+\omega_i\atop\mbox{\scriptsize (mod $Q$)}} 
\chi^{(i)}_{l,k}=\dfrac{(x^{2n};x^{2n})_\infty}{(x^2;x^2)_\infty} 
\left( \dfrac{(x^{2r};x^{2r})_\infty}{
(x^{2r-2};x^{2r-2})_\infty} \right)^{(n-1)(n-2)/2} G'_\xi, 
\label{eq:sumformula}
\end{equation}
where 
$$
G'_\xi =\prod_{\mu <\nu} [\xi_{\mu\nu}]'. 
$$
Eqs. (\ref{eq:sumformula}) and (\ref{eq:chi^i}) imply that 
\begin{equation}
\chi^{(i)}=\dfrac{1}{b_l} 
\sum_{k\equiv l+\omega_i\atop\mbox{\scriptsize (mod $Q$)}} 
\chi^{(i)}_{l,k}. 
\label{eq:chi-rel}
\end{equation}
where 
\begin{equation}
b_l =\left( \dfrac{
(x^{2r};x^{2r})_\infty}{(x^{2r-2};x^{2r-2})_\infty} \right)^{(n-1)(n-2)/2} 
G'_\xi . 
\label{eq:b_l}
\end{equation}

\subsection{Tail operator}

Let us introduce the dual intertwining vectors (see figure 2) 
satisfying 
\begin{equation}
\sum_{\mu =0}^{n-1} t_\mu^*  (v)^{a'}_{a}
t^\mu (v)^{a}_{a''} =\delta_{a''}^{a'}, ~~~~ 
\sum_{\nu =0}^{n-1} t^\mu (v)^{a}_{a-\bar{\varepsilon}_\nu} 
t_{\mu'}^* (v)^{a-\bar{\varepsilon}_\nu}_{a} =
\delta^\mu_{\mu'}. \label{eq:dual-t}
\end{equation}

\unitlength 1mm
\begin{picture}(100,20)
\put(40,3){\begin{picture}(101,0)
\put(-10,4){$\displaystyle\sum_{\mu=0}^{n-1}$}
\put(10.2,-3.){$a'$}
\put(-2.7,-2.8){$a$}
\put(4.15,0.1){\scriptsize{$\wedge$}}
\put(10,0){\vector(-1,0){10}}
\multiput(0,0)(0,2.2){5}{\line(0,1){1.2}}
\put(10.2,10.5){$a''$}
\put(-2.7,10.5){$a$}
\put(10,10){\vector(-1,0){10}}
\put(6,5){$\mu$}
\put(4.15,8.6){\scriptsize{$\vee$}}
\put(5,1.2){\line(0,1){7.6}}
\put(5,0){\line(0,-1){1}}
\put(5,-1.5){\vector(0,-1){2}}
\put(4.22,-6){$v$}
\put(15.,4.){$=\delta^{a''}_{a'}$,}
\end{picture}
}
\put(88,3){\begin{picture}(101,0)
\put(-10,4){$\displaystyle\sum_{a'}$}
\put(-3,4){$a$}
\put(11,4){$a'$}
\put(10,5){\vector(-1,0){10}}
\put(4.15,5.1){\scriptsize{$\wedge$}}
\put(4.15,3.6){\scriptsize{$\vee$}}
\put(5,6.4){\line(0,1){3.6}}
\put(6,7){$\mu'$}
\put(5,3.6){\vector(0,-1){3.6}}
\put(6,2){$\mu$}
\put(4.2,-2.2){$v$}
\put(17,4.){$=\delta^{\mu'}_{\mu}$.}
\end{picture}
}
\end{picture}

\vspace{3mm}

\begin{center}
Figure 2. Picture representation of the dual intertwining 
vectors. 
\end{center}

{}From (\ref{eq:Rtt=Wtt}) and (\ref{eq:dual-t}), we have (cf. figure 3) 
\begin{equation}
t^*(v_{1})^{b}_{c}\otimes t^*(v_{2})^{a}_{b}
R(v_{1}-v_2 )=
\displaystyle\sum_{d} 
W\left[ \left. \begin{array}{cc} 
c & d \\ b & a \end{array} \right| v_{1}-v_2 \right]
t^*(v_{1} )^{a}_{d}\otimes t^*(v_{2} )^{d}_{c}. 
\label{eq:dJMO}
\end{equation}

\unitlength 1mm
\begin{picture}(100,20)
\put(23,0){
\begin{picture}(101,0)
\put(20,3){\begin{picture}(101,0)
\put(10,0){\vector(-1,0){10}}
\put(0,0){\vector(0,1){10}}
\put(5,6.4){\line(0,1){3.6}}
\put(-.2,4.4){\scriptsize{$>$}}
\put(1.4,5){\line(1,0){10}}
\put(0,5){\line(-1,0){1}}
\put(-1.5,5){\vector(-1,0){2}}
\put(-6.5,4.2){$v_1$}
\put(4.15,.1){\scriptsize{$\wedge$}}
\put(5,1.4){\line(0,1){10}}
\put(5,0){\line(0,-1){1}}
\put(5,-1.5){\vector(0,-1){2}}
\put(4,-5.8){$v_2$}
\put(-2.5,10.5){$c$}
\put(10.5,-1.5){$a$}
\put(-2.5,-1.8){$b$}
\put(17,4){$=\;\displaystyle\sum_{d}$} 
\end{picture}
}
\put(56,3){\begin{picture}(101,0)
\put(10,0){\vector(-1,0){10}}
\put(10,0){\vector(0,1){10}}
\put(0,0){\vector(0,1){10}}
\put(10,10){\vector(-1,0){10}}
\put(9.8,4.4){\scriptsize{$>$}}
\multiput(10,5)(-2.2,0){6}{\line(-1,0){1.2}}
\put(11.6,5){\line(1,0){2}}
\put(-1.2,5){\vector(-1,0){2}}
\put(-7,4.2){$v_1$}
\put(4.15,10.1){\scriptsize{$\wedge$}}
\put(5,11.5){\line(0,1){2}}
\multiput(5,10)(0,-2.2){6}{\line(0,-1){1.2}}
\put(5,-1.2){\vector(0,-1){2}}
\put(4,-5.3){$v_2$}
\put(10.5,10.1){$d$}
\put(-2.5,10.5){$c$}
\put(10.5,-1.5){$a$}
\put(-2.5,-1.8){$b$}
\end{picture}
}
\end{picture}
}
\end{picture}

\vspace{2mm}

\begin{center}
Figure 3. Vertex-face correspondence by dual intertwining 
vectors. 
\end{center}

Now introduce the intertwining operators between 
${\cal H}^{(i)}$ and 
${\cal H}^{(i)}_{l,k}$ ($k=l+\omega_{i}$ (mod $Q$)): 
\begin{equation}
\begin{array}{rcl}
T(u ){}^{\xi a_0}&=&
\displaystyle\prod_{j=0}^\infty 
t^{\mu_j}(-u ){}^{a_j}_{a_{j+1}}: 
{\cal H}^{(i)}\rightarrow {\cal H}^{(i)}_{l,k}, \\
T(u ){}_{\xi a_0}&=&
\displaystyle\prod_{j=0}^\infty 
t^*_{\mu_j}(-u ){}_{a_j}^{a_{j+1}}: 
{\cal H}^{(i)}_{l,k}\rightarrow {\cal H}^{(i)}, 
\end{array}
\label{eq:T^_}
\end{equation}
where $k=a_0+\rho$ and $l=\xi +\rho$, and 
$0<\Re (u)<\frac{n}{2}+1$. 
Tail operator $\Lambda$ (see figure 4) is defined by 
\begin{equation}
\Lambda (u )_a^{a'}=T(u )^{\xi a'}T(u)_{\xi a}. 
\label{eq:L=TT}
\end{equation}
Let 
\begin{equation}
L\left[  \left. \begin{array}{cc} a'_0 & a'_1 \\
a_0 & a_1 \end{array} \right| u \right] :=
\sum_{\mu =0}^{n-1} t^*_\mu (-u)_{a_0}^{a_1} 
t^\mu (-u)^{a'_0}_{a'_1}. 
\label{eq:Lop}
\end{equation}
Then we have 
\begin{equation}
\Lambda(u ){}_{a_0}^{a'_0}=
\prod_{j=0}^\infty L\left[  \left. \begin{array}{cc} 
a'_j & a'_{j+1} \\
a_j & a_{j+1} \end{array} \right| u \right]. 
\label{eq:Lambda}
\end{equation}

\vspace{5mm}

\unitlength 1.4mm
\begin{picture}(100,20)
\put(-18,0){\begin{picture}(101,0)
\put(18,10){$\Lambda(u ){}_{a_0}^{a'_0}=$}
\multiput(65,5)(-10,0){3}{\vector(-1,0){10}}
\multiput(66,5)(2,0){5}{\line(1,0){1}}
\multiput(65,15)(-10,0){3}{\vector(-1,0){10}}
\multiput(66,15)(2,0){5}{\line(1,0){1}}
\put(110,5){\vector(-1,0){10}}
\put(110,15){\vector(-1,0){10}}
\put(100,5){\vector(-1,0){10}}
\put(100,15){\vector(-1,0){10}}
\put(90,5){\line(-1,0){4}}
\put(90,15){\line(-1,0){4}}
\put(86,5){\vector(-1,0){10}}
\put(86,15){\vector(-1,0){10}}
\put(34,2.3){$a_0$}
\put(44,2.3){$a_1$}
\put(54,2.3){$a_2$}
\put(64,2.3){$a_3$}
\put(34,16.){$a'_0$}
\put(44,16.2){$a'_1$}
\put(54,16.2){$a'_2$}
\put(64,16.2){$a'_3$}
\put(75,2.3){$\xi$}
\put(84,2.3){$\cdots$}
\put(89,2.3){$\xi\!\!+\!\!\omega_2$}
\put(99,2.3){$\xi\!\!+\!\!\omega_1$}
\put(109,2.3){$\xi$}
\put(75,16.){$\xi$}
\put(84,16){$\cdots$}
\put(89,16.){$\xi\!\!+\!\!\omega_2$}
\put(99,16){$\xi\!\!+\!\!\omega_1$}
\put(109,16.){$\xi$}
\put(39.4165,5.){\scriptsize{$\wedge$}}
\put(39.4165,13.9){\scriptsize{$\vee$}}
\put(40,6.0){\line(0,1){7.9}}
\put(40,5){\line(0,-1){1}}
\put(40,3.5){\vector(0,-1){2}}
\put(38.,-.3){$-u$}
\put(49.4165,5.){\scriptsize{$\wedge$}}
\put(49.4165,13.9){\scriptsize{$\vee$}}
\put(50,6.){\line(0,1){7.9}}
\put(59.4165,5.){\scriptsize{$\wedge$}}
\put(59.4165,13.9){\scriptsize{$\vee$}}
\put(60,6.){\line(0,1){7.9}}
\put(80.4165,5.){\scriptsize{$\wedge$}}
\put(80.4165,13.9){\scriptsize{$\vee$}}
\put(81,6.){\line(0,1){7.9}}
\put(94.4165,5.){\scriptsize{$\wedge$}}
\put(94.4165,13.9){\scriptsize{$\vee$}}
\put(95,6.){\line(0,1){7.9}}
\put(104.4165,5.){\scriptsize{$\wedge$}}
\put(104.4165,13.9){\scriptsize{$\vee$}}
\put(105,6.){\line(0,1){7.9}}
\end{picture}
}
\end{picture}

\vspace{3mm}

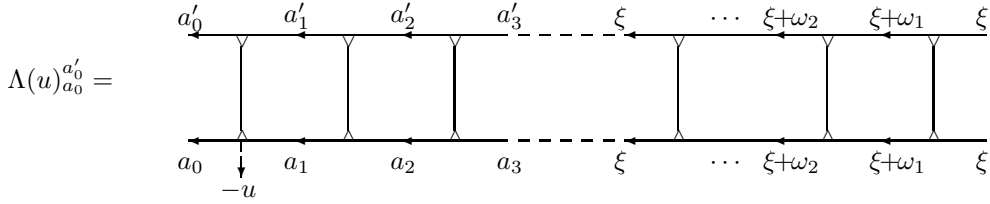
Figure 4. Tail operator $\Lambda(u ){}_{a_0}^{a'_0}$. 
The upper (resp. lower) half stands for $T(u ){}^{\xi a_0}$ (resp. 
$T(u ){}_{\xi a_0}$). 

\vspace{7mm}

Here we notice that in the `low temperature' 
limit $t^*_j (-u)_{\xi +\omega_{j+1}}^{\xi +\omega_{j}} 
t^j (-u)^{\xi +\omega_{j+1}}_{\xi +\omega_{j}}$ is much 
greater than other 
$t^*_\mu (-u)_{\xi +\omega_{j+1}}^{\xi +\omega_{j}} 
t^\mu (-u)^{\xi +\omega_{j+1}}_{\xi +\omega_{j}}$ ($\mu \neq j$). 

Note that 
\begin{equation}
L\left[ \left. \begin{array}{cc} a' & a'-\bar{\varepsilon}_\nu \\
a & a-\bar{\varepsilon}_\mu \end{array} \right| u \right] = 
\dfrac{[u+\bar{a}_\mu -\bar{a'}_\nu ]}{[u]}
\prod_{j\neq\mu} \dfrac{[\bar{a'}_\nu -\bar{a}_{j}]}
{[a_{\mu j}]}. 
\label{eq:Lop-ex}
\end{equation}
It is obvious from (\ref{eq:dual-t}), we have 
\begin{equation}
L\left[ \left. \begin{array}{cc} a & a' \\
a & a'' \end{array} \right| u \right] = \delta_{a''}^{a'}.  
\label{eq:L-inv}
\end{equation}
We therefore have 
\begin{equation}
\Lambda (u)^a_a=1. \label{eq:Lambda=1}
\end{equation}
{}From (\ref{eq:chi-rel}) and (\ref{eq:Lambda=1}), 
we may assume that 
\begin{equation}
\rho^{(i)}=\dfrac{1}{
b_l} \sum_{k\equiv l+\omega_i\atop\mbox{\scriptsize (mod $Q$)}} 
T (u)_{\xi a} \rho^{(i)}_{l,k} T(u)^{\xi a}. 
\label{eq:rho-rel}
\end{equation}

\subsection{Commutation relations among $\Lambda$ and $\phi$}

By using the vertex-face correspondence (see figure 5), we obtain 
\begin{equation}
T(u)^{\xi b} \Phi^\mu (v) =\sum_a t^\mu 
(v-u ){}_{a}^{b}\Phi (v)^b_aT(u)^{\xi a}, 
\label{eq:T^Phi}
\end{equation}
\begin{equation}
T(u)_{\xi b} \Phi (v)^b_a =\sum_\mu t^*_\mu 
(v-u ){}^{a}_{b}\Phi^\mu (v)T(u)_{\xi a}. 
\label{eq:T_phi}
\end{equation}

\unitlength 1mm
\begin{picture}(100,30)
\put(53,-42){\begin{picture}(101,0)
\put(40,55){\vector(0,1){10}}
\put(50,55){\vector(0,1){10}}
\put(50,55){\vector(-1,0){10}}
\put(50,65){\vector(-1,0){10}}
\put(60,55){\vector(0,1){10}}
\put(60,55){\vector(-1,0){10}}
\put(60,65){\vector(-1,0){10}}
\put(44.1,64.9){\scriptsize{$\wedge$}}
\put(54.1,64.9){\scriptsize{$\wedge$}}
\put(39,52){$a$}
\put(39,66){$b$}
\multiput(38,60)(2,0){13}{\line(1,0){1}}
\put(37.5,60){\vector(-1,0){2}}
\multiput(45,54)(0,2){6}{\line(0,1){1}}
\multiput(55,54)(0,2){6}{\line(0,1){1}}
\put(45,53.5){\vector(0,-1){1.5}}
\put(55,53.5){\vector(0,-1){1.5}}
\put(45,66){\line(0,1){3}}
\put(55,66){\line(0,1){3}}
\put(41.5,49.5){$-u$}
\put(51.5,49.5){$-u$}
\put(26.5,59.5){$v\!-\!u$}
\put(66,59){$=\;\displaystyle\sum_\mu$}
\end{picture}
}
\put(53,-42){\begin{picture}(101,0)
\put(80,55){\vector(0,1){10}}
\put(79,52){$a$}
\put(79,66){$b$}
\put(100,55){\vector(-1,0){10}}
\put(90,55){\vector(-1,0){10}}
\put(81.5,60){\line(1,0){18.5}}
\put(80.,59.4){\scriptsize{$>$}}
\thinlines
\put(82.,61){$\mu$}
\put(84.1,54.9){\scriptsize{$\wedge$}}
\put(85,56.3){\line(0,1){8.8}}
\multiput(85,54)(0,2){1}{\line(0,1){1}}
\put(85,53.5){\vector(0,-1){1.5}}
\put(81.5,49.5){$-u$}
\put(94.1,54.9){\scriptsize{$\wedge$}}
\put(95,56.3){\line(0,1){8.8}}
\multiput(95,54)(0,2){1}{\line(0,1){1}}
\put(95,53.5){\vector(0,-1){1.5}}
\put(91.5,49.5){$-u$}
\end{picture}
}
\put(5,-42){\begin{picture}(101,0)
\put(23,59){$=\;\displaystyle\sum_{a}$}
\put(40,55){\vector(0,1){10}}
\put(50,55){\vector(0,1){10}}
\put(50,55){\vector(-1,0){10}}
\put(50,65){\vector(-1,0){10}}
\put(60,55){\vector(0,1){10}}
\put(60,55){\vector(-1,0){10}}
\put(60,65){\vector(-1,0){10}}
\put(44.1,53.5){\scriptsize{$\vee$}}
\put(45,53.5){\vector(0,-1){2.5}}
\put(54.1,53.5){\scriptsize{$\vee$}}
\put(55,53.5){\vector(0,-1){2.5}}
\put(39,52){$a$}
\put(39,66){$b$}
\multiput(40,60)(2,0){12}{\line(1,0){1}}
\put(38.,59.4){\scriptsize{$<$}}
\put(38.3,60){\vector(-1,0){3.5}}
\put(36,61.5){$\mu$}
\multiput(45,55)(0,2){7}{\line(0,1){1}}
\multiput(55,55)(0,2){7}{\line(0,1){1}}
\put(41.5,49){$-u$}
\put(51.5,49){$-u$}
\end{picture}
}
\put(-70,-42){\begin{picture}(101,0)
\put(65.5,59){$v\!-\!u$}
\put(95,65){\vector(-1,0){10}}
\put(85,65){\vector(-1,0){10}}
\put(95,60){\vector(-1,0){21}}
\put(79.1,63.5){\scriptsize{$\vee$}}
\put(80,63.8){\vector(0,-1){8.8}}
\put(76.5,52.5){$-u$}
\put(89.1,63.5){\scriptsize{$\vee$}}
\put(90,63.8){\vector(0,-1){8.8}}
\put(86.5,52.5){$-u$}
\put(76.5,61){$\mu$}
\put(94,66){$k'$}
\end{picture}
}
\end{picture}

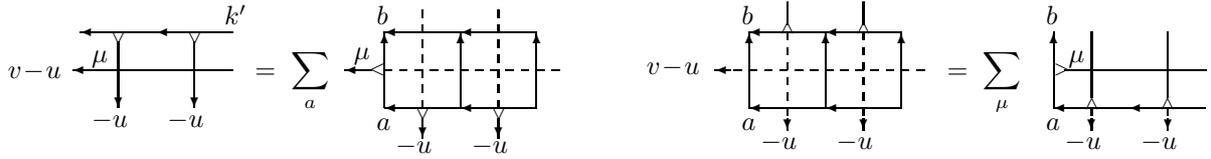
Figure 5. Commutation relations among $T(v_0)^{a\xi}$, 
$T(v_0)_{a\xi}$, and the type I vertex operators 
in vertex and face models. 

From these commutation relations and the definition of 
the tail operator (\ref{eq:L=TT}) we have 
\begin{equation}
\Lambda (u )^c_b\Phi (v)^b_a=
\sum_{d}L\left[ \left. \begin{array}{cc} c & d \\
b & a \end{array} \right| u-v \right]
\Phi (v)^c_d\Lambda (u)^d_a. 
\label{eq:Lambda-phi}
\end{equation}

\section{Vertex operator approach}

One of the most standard ways to calculate correlation functions 
is the vertex operator approach \cite{JMbk} 
on the basis of free field representation. 
In subsection 4.2 we recall the free field representation for 
the $A^{(1)}_{n-1}$ model \cite{AJMP}. The type I vertex operators 
of the $A^{(1)}_{n-1}$ model can be constructed in terms of 
basic bosons introduced in \cite{FL,AKOS}. 
The $A^{(1)}_{n-1}$ model has the so-called $\sigma$-invariance. 
The free field representation of type I vertex operator given in 
subsection 4.2 is not invariant under $\sigma$-transformation. Thus, 
we give other free field representations in subsection 4.3. 
We also need the bosonized CTM Hamiltonian 
of the $A^{(1)}_{n-1}$ model \cite{FHSW} in order to obtain 
correlation functions of the $A^{(1)}_{n-1}$ model. 
In subsection 4.4 we discuss the space of states of 
unrestricted $A^{(1)}_{n-1}$ model. The free field representation 
of the tail operator is presented in subsection 4.5. 

\subsection{Bosons}

Let us consider the bosons
$B_m^j\,(1\leqslant j \leqslant n-1, m \in \mathbb{Z}
\backslash \{0\})$
with the commutation relations
\begin{equation}
[B_m^j,B_{m'}^k]
=\left\{ \begin{array}{ll} 
m\dfrac{[(n-1)m]_x}{[nm]_x}
\dfrac{[(r-1)m]_x}{[rm]_x}\delta_{m+m',0}, & (j=k)\\
-mx^{{\rm sgn}(j-k)nm}\dfrac{[m]_x}{[nm]_x}
\dfrac{[(r-1)m]_x}{[rm]_x}\delta_{m+m',0}, & (j\neq k), 
\end{array} \right. 
\label{eq:comm-B}
\end{equation}
where the symbol $[a]_x$ stands for
$(x^a-x^{-a})/(x-x^{-1})$.
Define $B_m^n$ by
\begin{eqnarray*}
\sum_{j=1}^n x^{-2jm}B_m^j=0.
\end{eqnarray*}
Then the commutation relations (\ref{eq:comm-B}) 
holds for all $1\leqslant j,k \leqslant n$.
These oscillators were introduced in \cite{FL,AKOS}. 

For $\alpha, \beta \in {\mathfrak h}^*:=
{\mathbb C}\omega_0 \oplus {\mathbb C}\omega_1 \oplus \cdots 
{\mathbb C}\omega_{n-1}$, let us define 
the zero mode operators $P_\alpha, Q_\beta$ 
with the commutation relations 
\begin{equation*}
[P_{\alpha},\sqrt{-1}Q_{\beta}]=\langle \alpha,\beta \rangle, 
~~~~ [P_{\alpha}, B_m^j ]=
[Q_{\beta}, B_m^j ]=0. 
\end{equation*}

~\\
We will deal with the bosonic Fock spaces 
${\cal{F}}_{l,k}, (l,k \in {\mathfrak h}^*)$
generated by $B_{-m}^j (m>0)$
over the vacuum vectors $|l,k\rangle$ :
\begin{eqnarray*}
{\cal{F}}_{l,k}=
\mathbb{C}[\{ B_{-1}^j, B_{-2}^j,\cdots \}_{
1\leqslant j \leqslant n}]|l,k\rangle,
\end{eqnarray*}
where
\begin{eqnarray*}
B_m^j|l,k\rangle&=&0 ~(m>0),\\
P_{\alpha}|l,k\rangle &=&\langle \alpha,
\beta_1 k+\beta_2 l \rangle
|l,k\rangle,\\
|l,k\rangle&=&\exp \left(\sqrt{-1}(\beta_1Q_k+
\beta_2Q_l)\right)|0,0\rangle, 
\end{eqnarray*}
where $\beta_1$ and $\beta_2$ are defined by (\ref{eq:beta_12}).

\subsection{Type I vertex operators}

Let us define the basic operators for $j=1,\cdots,n-1$
\begin{eqnarray}
U_{-\alpha_j}(v)&=&z^{\frac{r-1}{r}}
:\exp\left(-\beta_1 (\sqrt{-1}Q_{\alpha_j}
+P_{\alpha_j}{\rm log} z)+\sum_{m \neq 0}\frac{1}{m}
(B_m^j-B_m^{j+1})(x^jz)^{-m}\right):, \\
U_{\omega_j}(v)&=&z^{\frac{r-1}{2r}\frac{j(n-j)}{n}}
:\exp\left(\beta_1 (\sqrt{-1}Q_{\omega_j}
+P_{\omega_j}{\rm log}z)
-\sum_{m\neq 0}\frac{1}{m} \sum_{k=1}^j 
x^{(j-2k+1)m}B_m^kz^{-m}\right):, 
\end{eqnarray}
where $\beta_1 =-\sqrt{\frac{r-1}{r}}$, and $z=x^{2v}$ 
as usual. 
Following commutation relations are useful: 
\begin{eqnarray}
U_{\omega_1}(v)U_{\omega_j}(v')&=&r_j (v-v')
U_{\omega_j}(v')U_{\omega_1}(v), \\
U_{-\alpha_j}(v)U_{\omega_j}(v')&=&-f(v-v',0)
U_{\omega_j}(v')U_{-\alpha_j}(v), \\
U_{-\alpha_j}(v)U_{-\alpha_{j+1}}(v')&=&-f(v-v',0)
U_{-\alpha_{j+1}}(v')U_{-\alpha_j}(v), \\
U_{-\alpha_j}(v)U_{-\alpha_j}(v')&=&g(v-v')
U_{-\alpha_j}(v')U_{-\alpha_j}(v). 
\end{eqnarray}

In the sequel we set
\begin{eqnarray*}
\pi_\mu=\sqrt{r(r-1)}P_{\bar{\varepsilon}_\mu},~
\pi_{\mu \nu}=\pi_\mu-\pi_\nu.
\end{eqnarray*}
The $\pi_{\mu \nu}$
acts on ${\cal{F}}_{l,k}$ as a scalor 
$\langle \varepsilon_\mu -\varepsilon_\nu,
rl-(r-1)k\rangle$.

For $0 \leqslant \mu \leqslant n-1$ 
define the type I vertex operator \cite{AJMP} by
\begin{equation}
\begin{array}{cl}
&\phi_\mu(v)=\displaystyle\oint
\prod_{j=1}^{\mu}\frac{dz_j}{2\pi \sqrt{-1} z_j}
U_{\omega_1}(v)U_{-\alpha_1}(v_1)\cdots 
U_{-\alpha_\mu}(v_{\mu})
\prod_{j=0}^{\mu-1}f(v_{j+1}-v_{j},\pi_{j \mu}) 
\prod_{j=0\atop j\neq\mu}^{n-1} [\pi_{j\mu }]^{-1} \\
=&\displaystyle (-1)^\mu \oint
\prod_{j=1}^{\mu}\frac{dz_j}{2\pi \sqrt{-1} z_j}
U_{-\alpha_\mu}(v_{\mu})\cdots U_{-\alpha_1}(v_1)
U_{\omega_1}(v)
\prod_{j=0}^{\mu-1}f(v_{j}-v_{j+1},1-\pi_{j \mu}) 
\prod_{j=0\atop j\neq\mu}^{n-1} [\pi_{j\mu }]^{-1}, 
\label{eq:type-I}
\end{array}
\end{equation}
where $v_0 =v$ and $z_j=x^{2v_j}$. 
The integral contour for $z_j$-integration 
encircles the poles at $z_j=x^{1+2kr}z_{j-1}\,(k 
\in \mathbb{Z}_{\geqslant 0})$, but not the poles at 
$z_j=x^{-1-2kr}z_{j-1}\,(k \in 
\mathbb{Z}_{\geqslant 0})$, for $1\leqslant j\leqslant \mu$. 

Note that 
\begin{equation}
\phi_\mu(v): {\cal{F}}_{l,k} 
\longrightarrow {\cal{F}}_{l,k+\bar{\varepsilon}_\mu}. 
\label{eq:k-shift}
\end{equation}
These type I vertex operators satisfy 
the following commutation relations on 
${\cal{F}}_{l,k}$: 
\begin{eqnarray}
\phi_{\mu_1}(v_1)\phi_{\mu_2}(v_2)
=\sum_{\varepsilon_{\mu_1}+\varepsilon_{\mu_2}
=\varepsilon_{\mu_1'}+\varepsilon_{\mu_2'} }W\left[\left.
\begin{array}{cc}
a+\bar{\varepsilon}_{\mu_1}+\bar{\varepsilon}_{\mu_2}&
a+\bar{\varepsilon}_{\mu_1'}\\
a+\bar{\varepsilon}_{\mu_2}&a
\end{array}\right|v_1-v_2 \right]
{\phi}_{\mu_2'}(v_2)
{\phi}_{\mu_1'}(v_1). 
\label{eq:CR-I}
\end{eqnarray}
We thus denote the operator $\phi_\mu (v)$ by 
$\Phi (v )^{a+\bar{\varepsilon}_\mu}_a$ 
on the bosonic Fock space ${\cal{F}}_{l,a+\rho}$. 
We notice that our vertex operator (\ref{eq:type-I}) has 
different normalization from that originally constructed 
in \cite{AJMP} because of the difference of the Boltzmann 
weight $W$. Furthermore, the range of $\mu$ is shifted from 
that of \cite{AJMP} by $1$ so that our $\phi_\mu (v)$ corresponds 
to $\phi_{\mu +1} (v)$ in \cite{AJMP}, up to normalization. 

Dual vertex operators are likewise defined as follows: 
\begin{equation}
\begin{array}{cl}
&\phi^*_\mu(v)=c_n^{-1}\displaystyle\oint
\prod_{j=\mu +1}^{n-1}\frac{dz_j}{2\pi \sqrt{-1} z_j}
U_{\omega_{n-1}}\left(v-\frac{n}{2} \right)U_{-\alpha_{n-1}}(v_{n-1})\cdots 
U_{-\alpha_{\mu +1}}(v_{\mu +1})
\prod_{j=\mu +1}^{n-1}f(v_{j}-v_{j+1},\pi_{\mu j}) \\
=&c_n^{-1}\displaystyle (-1)^{n-1-\mu} \oint
\prod_{j=\mu +1}^{n-1}\frac{dz_j}{2\pi \sqrt{-1} z_j}
U_{-\alpha_{\mu +1}}(v_{\mu +1})\cdots U_{-\alpha_{n-1}}(v_{n-1})
U_{\omega_{n-1}}\left(v-\frac{n}{2}\right)
\prod_{j=\mu +1}^{n-1}f(v_{j+1}-v_{j},1-\pi_{\mu j}) 
\label{eq:type-I*}
\end{array}
\end{equation}
where $v_n =v-\tfrac{n}{2}$, and 
$$
c_n =x^{\frac{r-1}{r}\frac{n-1}{2n}} \dfrac{g_{n-1}(x^n)}{
(x^2; x^{2r})_\infty^n (x^{2r},x^{2r})_\infty^{2n-3}}. 
$$ 
The integral contour for $z_j$-integration 
encircles the poles at $z_j=x^{1+2kr}z_{j+1}\,(k 
\in \mathbb{Z}_{\geqslant 0})$, but not the poles at 
$z_j=x^{-1-2kr}z_{j+1}\,(k \in 
\mathbb{Z}_{\geqslant 0})$, for $\mu +1\leqslant j\leqslant n-1$. 
Note that 
\begin{equation}
\phi^*_\mu (v): {\cal{F}}_{l,k} 
\longrightarrow {\cal{F}}_{l,k-\bar{\varepsilon}_\mu}. 
\label{eq:k-shift*}
\end{equation}
The operators $\phi_\mu (v)$ and $\phi^*_\mu (v)$ are dual 
in the following sense: 
\begin{equation}
\sum_{\mu =0}^{n-1} \phi^*_\mu (v)\phi_\mu (v)=1. 
\label{eq:dualrel}
\end{equation}
We notice that our dual vertex operator $\phi^*_\mu (v)$ 
coincides with 
$\bar{\phi}^{^*(n-1)}_{\mu +1}\left(v-\frac{n}{2}\right)$ in 
\cite{AJMP}. 

\subsection{Other representations} 

The present face model has the so-called $\sigma$-invariance: 
$$
W\left[\left.
\begin{array}{cc}
\sigma (c)&\sigma (d) \\
\sigma (b)&\sigma (a)
\end{array}\right|v \right]=W\left[\left.
\begin{array}{cc}
c&d \\
b&a
\end{array}\right|v \right], ~~~~ \sigma (\omega_\mu )=
\omega_{\mu +1}. 
$$
The free field representation (\ref{eq:type-I}) 
is not invariant under $\sigma$-transformation, so that 
we have other free field representations: 
\begin{equation}
\begin{array}{cl}
&\phi_{i+\mu }(v)=\displaystyle\oint
\prod_{j=1}^{\mu}\frac{dz_j}{2\pi \sqrt{-1} z_j}
U_{\omega_1}(v)U_{-\alpha_1}(v_1)\cdots 
U_{-\alpha_\mu}(v_{\mu})
\prod_{j=0}^{\mu-1}f(v_{j+1}-v_{j},\pi_{i+j i+\mu}) 
\prod_{j=0\atop j\neq\mu}^{n-1} [\pi_{i+j i+\mu }]^{-1} \\
=&\displaystyle (-1)^\mu \oint
\prod_{j=1}^{\mu}\frac{dz_j}{2\pi \sqrt{-1} z_j}
U_{-\alpha_\mu}(v_{\mu})\cdots U_{-\alpha_1}(v_1)
U_{\omega_1}(v)
\prod_{j=0}^{\mu-1}f(v_{j}-v_{j+1},1-\pi_{i+j i+\mu}) 
\prod_{j=0\atop j\neq\mu}^{n-1} [\pi_{i+j i+\mu }]^{-1}, 
\label{eq:type-I(i)}
\end{array}
\end{equation}
where $v_0 =v$ and $z_j=x^{2v_j}$, and the 
integral contours are same one as (\ref{eq:type-I}). 
In this representation the space of states 
${\cal H}^{(i)}_{l,k}$ should be identified with 
${\cal F}_{\sigma^{-i}(l),\sigma^{-i}(k)}$. 

\subsection{Free field realization of CTM Hamiltonian}

Let 
\begin{equation}
\begin{array}{rcl}
H_F &=&\displaystyle\sum_{m=1}^\infty 
\dfrac{[rm]_x}{[(r-1)m]_x} 
\sum_{j=1}^{n-1}\sum_{k=1}^j x^{(2k-2j-1)m} 
B_{-m}^k (B_m^j -B_m^{j+1}) +\dfrac{1}{2}
\sum_{j=1}^{n-1} P_{\omega_j}P_{\alpha_j} \\
&=&\displaystyle\sum_{m=1}^\infty 
\dfrac{[rm]_x}{[(r-1)m]_x} 
\sum_{j=1}^{n-1}\sum_{k=1}^j x^{(2j-2k-1)m} 
(B_{-m}^j -B_{-m}^{j+1}) B_{m}^k +\dfrac{1}{2}
\sum_{j=1}^{n-1} P_{\omega_j}P_{\alpha_j}
\end{array}
\label{eq:CTM-Fock}
\end{equation}
be the CTM Hamiltonian on the Fock space ${\cal F}_{l,k}$ \cite{FHSW}. 
Then we have the homogeneity relation
\begin{equation}
\phi_\mu (z) q^{H_F} =q^{H_F}\phi_\mu (q^{-1}z)
\label{eq:homo}
\end{equation}
and 
\begin{equation}
\mbox{tr}_{{\cal F}_{l,k}}\, \left( x^{2n H_F} G_a  \right) 
=\dfrac{x^{n |\beta_1 k+\beta_2l|^2}}
{(x^{2n};x^{2n})^{n-1}_\infty}G_a. 
\label{eq:chi-Fock}
\end{equation}
By comparing (\ref{eq:chi_lk}) and (\ref{eq:chi-Fock}), 
we conclude that $\rho^{(i)}_{l,k}=G_a x^{2nH_F}$ and 
${\cal H}^{(i)}_{l,k}={\cal F}_{l,k}$, where 
$k=a+\rho$. 

The relation between $\rho^{(i)}$ and $\rho^{(i)}_{l,k}$ is as follows: 
\begin{equation}
\rho^{(i)}= \sum_{
k\equiv l+\omega_i\atop\mbox{\scriptsize (mod $Q$)}} 
T(u)_{\xi a} \dfrac{\rho^{(i)}_{l,k}}{b_l} T(u)^{\xi a}. 
\label{eq:rel-rho^2}
\end{equation}

\subsection{Free field realization of tail operators}

Consider (\ref{eq:Lambda-phi}) for $(c,b,a)\rightarrow 
(a,a+\bar{\varepsilon}_{0}+\bar{\varepsilon}_{\mu}, 
a-\bar{\varepsilon}_{\mu})$, where $\mu\neq 0$. 
The coefficient $L$ diverges when $u\rightarrow v$, so that 
we obtain the following necessary condition: 
\begin{equation}
\prod_{j=1\atop j\neq\mu}^{n-1} [a_{0j}] 
\Phi (v)^a_{a-\bar{\varepsilon}_{0}} 
\Lambda (v)^{a-\bar{\varepsilon}_{0}}_{a-\bar{\varepsilon}_{\mu}}+
\prod_{j=1\atop j\neq\mu}^{n-1} 
[a_{\mu j}] \Phi (v)^a_{a-\bar{\varepsilon}_{\mu}} 
\Lambda (v)^{a-\bar{\varepsilon}_{\mu}}_{a-\bar{\varepsilon}_{\mu}}
=0. \label{eq:phiLambda}
\end{equation}
By solving (\ref{eq:phiLambda}), we obtain 
\begin{equation}
\begin{array}{rcl}
\Lambda (u)_{a-\bar{\varepsilon}_{\mu}}^{a-\bar{\varepsilon}_{0}}
&=& G_\pi
\displaystyle\oint \prod_{j=1}^{\mu} 
\dfrac{dz_j}{2\pi\sqrt{-1}z_j} 
U_{-\alpha_{1}}(v_{1}) \cdots 
U_{-\alpha_{\mu}}(v_{\mu}) \displaystyle\prod_{j=0}^{\mu -1} 
f(v_{j+1}-v_j, \pi_{j\mu }) G_\pi^{-1}, 
\end{array}
\label{eq:Bose-Lambda}
\end{equation}
where 
$$
G_\pi :=\displaystyle\prod_{\kappa <\lambda}[\pi_{\kappa\lambda}]. 
$$
Note that a free field representation of 
$\Lambda (u)_{a-\bar{\varepsilon}_{\mu}}^{a-\bar{\varepsilon}_{\nu}}$ 
for $\nu >0$ can be constructed on 
${\cal F}_{\sigma^{-\nu}(l),\sigma^{-\nu}(k)}$. 

In the next section we need a tail operator 
$\Lambda (u)_{a-\sum_{j=1}^N \bar{\varepsilon}_{\mu_j}}^
{a-\sum_{j=1}^N \bar{\varepsilon}_{\nu_j}}$ in order to calculate 
$n$-point functions. This type tail operator can be represented 
in terms of free bosons. In order to show this fact, let us 
introduce the symbol $\lesssim$ as follows. We say 
$\mu \lesssim \nu$ if $0\leqq \mu_0 \leqq \nu_0 \leqq n-1$ and 
$\mu = \mu_0$ (mod $n$), $\nu = \nu_0$ (mod $n$). 

It is clear that there exists $0\leqq i\leqq n-1$ 
such that 
$$
\sharp\{ j| \nu_j +i \lesssim 0 \} >0, 
$$
and 
$$
\sharp\{ j| \mu_j +i \lesssim m \} \leqq \sharp\{ j| \nu_j +i \lesssim m \}, 
$$
for every $0\leqq m\leqq n-1$. 
In this case a free field representation of the tail operator 
$\Lambda (u)_{a-\sum_{j=1}^N \bar{\varepsilon}_{\mu_j}}^
{a-\sum_{j=1}^N \bar{\varepsilon}_{\nu_j}}$ can be constructed on 
${\cal F}_{\sigma^{-i}(l),\sigma^{-i}(k)}$. 

\section{Correlation functions}

\subsection{General formulae}

Consider the local state probability (LSP) such that the state variable 
at $j$th site is equal to $\mu_j$ ($1\leqq j\leqq N$), under a certain 
fixed boundary condition. 
In order to obtain LSP, it is convenient to 
divide the lattice into four transfer matrices and $2N$ vertex operators 
as follows: 

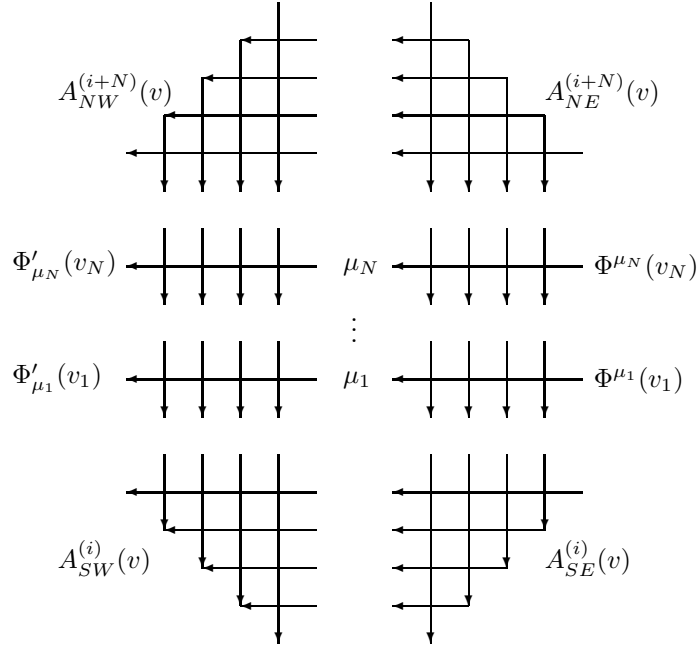
\unitlength 0.5mm
\begin{picture}(300,190)
\put(30,0){\begin{picture}(100,0)
\put(120,0){\begin{picture}(100,0)
\put(20,10){\vector(-1,0){20}}
\put(30,20){\vector(-1,0){30}}
\put(40,30){\vector(-1,0){40}}
\put(50,40){\vector(-1,0){50}}
\put(10,50){\vector(0,-1){50}}
\put(20,50){\vector(0,-1){40}}
\put(30,50){\vector(0,-1){30}}
\put(40,50){\vector(0,-1){20}}
\put(40,20){$A^{(i)}_{SE}(v)$}
\end{picture}
}
\put(120,120){\begin{picture}(100,0)
\put(20,40){\vector(-1,0){20}}
\put(30,30){\vector(-1,0){30}}
\put(40,20){\vector(-1,0){40}}
\put(50,10){\vector(-1,0){50}}
\put(10,50){\vector(0,-1){50}}
\put(20,40){\vector(0,-1){40}}
\put(30,30){\vector(0,-1){30}}
\put(40,20){\vector(0,-1){20}}
\put(40,24){$A^{(i+N)}_{NE}(v)$}
\end{picture}
}
\put(50,120){\begin{picture}(100,0)
\put(50,40){\vector(-1,0){20}}
\put(50,30){\vector(-1,0){30}}
\put(50,20){\vector(-1,0){40}}
\put(50,10){\vector(-1,0){50}}
\put(40,50){\vector(0,-1){50}}
\put(30,40){\vector(0,-1){40}}
\put(20,30){\vector(0,-1){30}}
\put(10,20){\vector(0,-1){20}}
\put(-18,24){$A^{(i+N)}_{NW}(v)$}
\end{picture}
}
\put(50,0){\begin{picture}(100,0)
\put(50,10){\vector(-1,0){20}}
\put(50,20){\vector(-1,0){30}}
\put(50,30){\vector(-1,0){40}}
\put(50,40){\vector(-1,0){50}}
\put(40,50){\vector(0,-1){50}}
\put(30,50){\vector(0,-1){40}}
\put(20,50){\vector(0,-1){30}}
\put(10,50){\vector(0,-1){20}}
\put(-18,20){$A^{(i)}_{SW}(v)$}
\end{picture}
}
\put(110,0){\begin{picture}(100,0)
\put(-10,70){\vector(-1,0){50}}
\put(-20,80){\vector(0,-1){20}}
\put(-30,80){\vector(0,-1){20}}
\put(-40,80){\vector(0,-1){20}}
\put(-50,80){\vector(0,-1){20}}
\put(-3,69){$\mu_1$}
\put(-1,80){$\vdots$}
\put(-90,69){$\Phi'_{\mu_1} (v_1)$}
\end{picture}
}
\put(110,30){\begin{picture}(100,0)
\put(-10,70){\vector(-1,0){50}}
\put(-20,80){\vector(0,-1){20}}
\put(-30,80){\vector(0,-1){20}}
\put(-40,80){\vector(0,-1){20}}
\put(-50,80){\vector(0,-1){20}}
\put(-3,69){$\mu_N$}
\put(-90,69){$\Phi'_{\mu_N} (v_N)$}
\end{picture}
}
%
\put(180,0){\begin{picture}(100,0)
\put(-10,70){\vector(-1,0){50}}
\put(-20,80){\vector(0,-1){20}}
\put(-30,80){\vector(0,-1){20}}
\put(-40,80){\vector(0,-1){20}}
\put(-50,80){\vector(0,-1){20}}
\put(-7,68){$\Phi^{\mu_1} (v_1)$}
\end{picture}
}
\put(180,30){\begin{picture}(100,0)
\put(-10,70){\vector(-1,0){50}}
\put(-20,80){\vector(0,-1){20}}
\put(-30,80){\vector(0,-1){20}}
\put(-40,80){\vector(0,-1){20}}
\put(-50,80){\vector(0,-1){20}}
\put(-7,68){$\Phi^{\mu_N} (v_N)$}
\end{picture}
}
\end{picture}
}
\end{picture}

\vspace{10mm}

\noindent Here, the incoming vertex operator $\Phi'_{\mu} (v)$ should be 
distinguished form the outgoing vertex operator 
$\Phi^{\mu} (v)$. 

Let us consider the normalized partition function 
with fixed $\mu_1, \cdots \mu_N$: 
\begin{equation}
\begin{array}{cl}
& P^{(i)}_{\mu_1 \cdots \mu_N}(v_1 , \cdots , v_N ) \\
:=& \displaystyle
\dfrac{1}{\chi^{(i)}}\mbox{tr}_{{\cal H}^{(i)}}\, 
\left( A^{(i)}_{SW}(v)\Phi'_{\mu_1} (v_1)\cdots 
\Phi'_{\mu_N} (v_N)A^{(i+N)}_{NW}(v)
A^{(i+N)}_{NE}(v)\Phi^{\mu_N} (v_N)\cdots 
\Phi^{\mu_1} (v_1)A^{(i)}_{SE}(v) \right). 
\end{array}
\end{equation}
In the vertex operator approach \cite{JMbk}, the LSP can be given by 
$P^{(i)}_{\mu_1 \cdots \mu_N}(v_1 , \cdots , v_N )|_{
v_1 =\cdots =v_N =v}$. In what follows, we denote 
$P^{(i)}_{\mu_1 \cdots \mu_N}=
P^{(i)}_{\mu_1 \cdots \mu_N}(v_1 , \cdots , v_N )|_{
v_1 =\cdots =v_N =v}$. 

The YBE and the crossing symmetry imply the following relation 
\cite{SPn}: 
\begin{equation}
\Phi'_{\mu} (v')A^{(i+1)}_{NW}(v)A^{(i+1)}_{NE}(v)
=A^{(i)}_{NW}(v)A^{(i)}_{NE}(v)\Phi^*_{\mu} (v'). 
\end{equation}

Thus, one point local state probability 
of $(\mathbb{Z}/n\mathbb{Z})$-symmetric model can be given by 
\begin{equation}
\begin{array}{rcl}
P_j^{(i)}&=&\dfrac{1}{\chi^{(i)}}\mbox{tr}_{{\cal H}^{(i)}}\, 
\left( \Phi^*_j (v) \Phi^j (v) \rho^{(i)} \right) \\
&=& \displaystyle
\dfrac{1}{\chi^{(i)}} \sum_{
k\equiv l+\omega_i\atop\mbox{\scriptsize (mod $Q$)}} 
\mbox{tr}_{{\cal H}^{(i)}_{l,k}}\, 
\left( T(u )^{\xi a}  
\Phi^*_j (v) \Phi^j (v) T(u)_{\xi a} 
\dfrac{\rho^{(i)}_{l,k}}{b_l} \right) \\
&=&\displaystyle
\dfrac{1}{\chi^{(i)}} \sum_{
k\equiv l+\omega_i\atop\mbox{\scriptsize (mod $Q$)}} 
\sum_{\mu ,\nu} t^*_j(v-u)^{
a}_{a+\bar{\varepsilon}_\nu} 
t^j(v-u)^{a+\bar{\varepsilon}_\nu}_{
a+\bar{\varepsilon}_\nu-\bar{\varepsilon}_\mu} \\ 
&& ~~~ \times \mbox{tr}_{{\cal H}^{(i)}_{l,k}}\, 
\left( \Phi^* (v)^{a}_{a+\bar{\varepsilon}_\nu} 
\Phi (v)^{a+\bar{\varepsilon}_\nu}_{
a+\bar{\varepsilon}_\nu-\bar{\varepsilon}_\mu } 
\Lambda (u)_{a}^{a+\bar{\varepsilon}_\nu-\bar{\varepsilon}_\mu
} \dfrac{\rho^{(i)}_{l,k}}{b_l} \right). 
\end{array}
\label{eq:Corr}
\end{equation}
Here in the third equality of (\ref{eq:Corr}), we use 
(\ref{eq:T^Phi}) and the fact 
that $\Phi^*_j (v)$, $\Phi^* (v)^{a}_{a+\bar{\varepsilon}_\nu}$ 
and $t^*_j(v)^{a}_{a+\bar{\varepsilon}_\nu}$ 
are given by the fusion of $n-1$ $\Phi^k (v)$'s, 
$\Phi (v)^{a+\bar{\varepsilon}_\mu}_{a}$'s and 
$t^k(v)^{a+\bar{\varepsilon}_\mu}_{a}$'s, respectively. 

In general, $N$-point local state possibility of this model 
can be given by 
\begin{equation}
\begin{array}{rcl}
P_{j_N\cdots j_1}^{(i)}&=&\dfrac{1}{\chi^{(i)}}\mbox{tr}_{{\cal H}^{(i)}}\, 
\left( \Phi^*_{j_N} (v)\cdots\Phi^*_{j_1} (v) \Phi^{j_1} (v)\cdots 
\Phi^{j_N} (v) \rho^{(i)} \right) \\
&=&\dfrac{1}{\chi^{(i)}} \displaystyle\sum_{
k\equiv l+\omega_i\atop\mbox{\scriptsize (mod $Q$)}} 
\sum_{a_1\cdots a_N\atop a'_1\cdots a'_N} t^*_{j_N}(v-u)^{a}_{a_N} \cdots t^*_{j_1}(v-u)^{a_2}_{a_1} 
t^{j_1}(v-u)^{a_1}_{a'_1} \cdots t^{j_N}(v-u)^{a'_{N-1}}_{a'_N} \\ 
&& ~~~ \times \mbox{tr}_{{\cal H}^{(i)}_{l,k}}\, 
\left( \Phi^* (v)^{a}_{a_N} \cdots \Phi^* (v)^{a_2}_{a_1} 
\Phi (v)^{a_1}_{a'_1} \cdots \Phi (v)^{a'_{N-1}}_{a'_N}
\Lambda (u)_{a}^{a'_N} \dfrac{\rho^{(i)}_{l,k}}{b_l} \right), 
\end{array}
\label{eq:Corr_N}
\end{equation}
where the second sum on the second line should be taken such that 
$(a,a_N)$, $\cdots$, $(a_2, a_1)$ and $(a'_1, a_1)$, $\cdots$, 
$(a'_N, a'_{N-1})$ are all admissible.

\subsection{Spontaneous polarization}

In this subsection we reproduce the expression for 
spontaneous polarization \cite{SPn}: 
\begin{equation}
\langle g \rangle ^{(i)} =\sum_{j=0}^{n-1} 
\omega^j P_j^{(i)}=
\omega ^{i+1} 
\displaystyle\frac{(x^2 ; x^{2})_{\infty}^{2}} 
                  {(x^{2r} ; x^{2r})_{\infty}^{2}} 
\displaystyle\frac{(\omega x^{2r} ; x^{2r})_{\infty}
                   (\omega ^{-1} x^{2r}; x^{2r})_{\infty}}
                  {(\omega x^2 ; x^{2})_{\infty}
                   (\omega ^{-1} x^{2}; x^{2})_{\infty}}, 
\label{eq:SPn}
\end{equation}
by performing traces and $n$-fold integrals on (\ref{eq:Corr}). 
In \cite{SPn} the expression (\ref{eq:SPn}) was obtained by 
solving a system of difference equations, the quantum 
Knizhnik--Zamolodchikov equations of level $-2n$. 

First we replace $a+\bar{\varepsilon}_\nu$ by $a$ for simplicity: 
\begin{equation}
\begin{array}{rcl}
P_j^{(i)}
&=&\displaystyle
\dfrac{1}{\chi^{(i)}} \sum_{
k\equiv l+\omega_{i+1}\atop\mbox{\scriptsize (mod $Q$)}} 
\sum_{\mu ,\nu} t^*_j(v-u)^{a-\bar{\varepsilon}_\nu}_{a} 
t^j(v-u)^{a}_{a-\bar{\varepsilon}_\mu} \\ 
&& ~~~ \times \mbox{tr}_{{\cal H}^{(i+1)}_{l,k}}\, 
\left( \Phi (v)_{a-\bar{\varepsilon}_\mu}^{a}
\Lambda (u)_{a-\bar{\varepsilon}_\nu}^{a-\bar{\varepsilon}_\mu} 
\dfrac{\rho^{(i)}_{l,k-\bar{\varepsilon}_\nu}}{b_l} 
\Phi^* (v)_{a}^{a-\bar{\varepsilon}_\nu} \right). 
\end{array}
\label{eq:Corr_1}
\end{equation}

We note that 
\begin{equation}
\sum_{j=0}^{n-1} \omega^j t^*_j(v-u)^{
a-\bar{\varepsilon}_{\nu}}_{a} 
t^j(v-u)^{a}_{
a-\bar{\varepsilon}_\mu }=
\dfrac{\mbox{[}v-u +a_{\mu \nu}\mbox{]}_\omega}
{\mbox{[}v-u\mbox{]}}\prod_{j=0\atop j\neq\nu}^{n-1} 
\dfrac{\mbox{[}a_{\mu j}\mbox{]}_\omega}{\mbox{[}a_{\nu j}\mbox{]}}, 
\label{eq:omega,t^*t}
\end{equation}
where 
$$
[v]_{\omega} 
=x^{\frac{v^2}{r}-v}\Theta_{x^{2r}} (\omega x^{2v}). 
$$

Thus, the spontaneous polarization can be reduced as 
$$
\langle g \rangle ^{(i)} =
\dfrac{1}{\chi^{(i)}} \sum_{\mu =0}^{n-1} \langle g \rangle ^{(i)}_\mu , 
$$
where 
$$
\langle g \rangle ^{(i)}_\mu =\sum_{
k\equiv l+\omega_{i+1}\atop\mbox{\scriptsize (mod $Q$)}} 
\sum_{\nu =0}^{n-1} \dfrac{\mbox{[}v-u +a_{\mu \nu}\mbox{]}_\omega}
{\mbox{[}v-u\mbox{]}}\prod_{j=0\atop j\neq\nu}^{n-1} 
\dfrac{\mbox{[}a_{\mu j}\mbox{]}_\omega}{\mbox{[}a_{\nu j}\mbox{]}}
\mbox{tr}_{{\cal H}^{(i+1)}_{l,k}}\, 
\left( \Phi (v)_{a-\bar{\varepsilon}_\mu}^{a}
\Lambda (u)_{a-\bar{\varepsilon}_\nu}^{a-\bar{\varepsilon}_\mu} 
G_{a-\bar{\varepsilon}_\nu} 
\Phi^* \left(v+n\right)_{a}^{a-\bar{\varepsilon}_\nu} 
\dfrac{x^{2nH^{(i+1)}_{l,k}}}{b_l} \right). 
$$

When $\mu =0$, 
in order to calculate the operator product 
$\Phi (v)_{a-\bar{\varepsilon}_0}^{a}
\Lambda (u)_{a-\bar{\varepsilon}_\nu}^{a-\bar{\varepsilon}_0
} G_{a-\bar{\varepsilon}_\nu} 
\Phi^* \left(v+n\right)_{a}^{a-\bar{\varepsilon}_\nu}$, 
the following operator product formulae 
are useful: 
\begin{equation}
\begin{array}{cl}
& \displaystyle 
c_n^{-1}U_{\omega_1}(v) U_{-\alpha_1}(v_1 )\cdots 
U_{-\alpha_{n-1}}(v_{n-1} )U_{\omega_{n-1}}\left(v+\tfrac{n}{2}\right) \\
=& x^{-\frac{n-1}{2}\frac{r-1}{r}} (x^2; x^{2r})_\infty^n 
(x^{2r}; x^{2r})_\infty^{2n-3}
z^{-\frac{r-1}{nr}} 
\displaystyle \prod_{j=0}^{n-1} z_{j}^{-\frac{r-1}{r}} 
\dfrac{(x^{2r-1}z_{j+1}/z_{j}; x^{2r})_\infty}{
(xz_{j+1}/z_{j}; x^{2r})_\infty} \\ 
\times& :U_{\omega_1}(v) U_{-\alpha_1}(v_1 )\cdots 
U_{-\alpha_{n-1}}(v_{n-1} )U_{\omega_{n-1}}\left(v+\frac{n}{2}\right):, 
\end{array}
\label{eq:U^{n+1}}
\end{equation}
where $z_0 =z$, and $z_n =x^nz$. Using (\ref{eq:U^{n+1}}) 
we have the following trace formulae: 
\begin{equation}
\begin{array}{cl}
& \mbox{tr}_{{\cal H}^{(i+1)}_{l,k}}\, 
\left( c_n^{-1}U_{\omega_1}(v) 
U_{-\alpha_1}(v_1 )\cdots U_{-\alpha_{n-1}}(v_{n-1} )
U_{\omega_{n-1}}\left(v+\tfrac{n}{2}\right)
x^{2nH^{(i+1)}_{l,k}} \right) \\
=& x^{n\left(\frac{r-1}{r}|k^2|-2\langle k,l\rangle 
+\frac{r}{r-1}|l|^2\right)} x^{2\frac{r-1}{r} 
\sum_{j=1}^{n-1} a_{0\,j}\left(v_{j+1}-v_j-\frac{1}{2}\right)
-2\sum_{j=1}^{n-1} \xi_{0\,j}\left(v_{j+1}-v_j-\frac{1}{2}\right)}\\
\times & (x^{2n}; x^{2n})_\infty 
(x^{2r}; x^{2r})_\infty^{2n-3} \dfrac{(x^2; x^{2n}, x^{2r})^n_\infty}{
(x^{2r+2n-2}; x^{2n}, x^{2r})^n_\infty} 
\\ 
\times & \displaystyle 
\prod_{j=0}^{n-1} 
\dfrac{(x^{2r-1}z_{j+1}/z_{j}; x^{2n}, x^{2r})_\infty 
(x^{2r+2n-1}z_{j}/z_{j+1}; x^{2n}, x^{2r})_\infty}{
(xz_{j+1}/z_{j}; x^{2n}, x^{2r})_\infty 
(x^{2n+1}z_{j}/z_{j+1}; x^{2n}, x^{2r})_\infty}. 
\end{array}
\label{eq:*trace}
\end{equation}

Let us denote the r.h.s. of (\ref{eq:*trace}) by 
$A^{(i)}_{l,k}(v; v_1, \cdots , v_{n-1})$. 
Then we have 
\begin{equation}
\begin{array}{rcl}
\langle g \rangle ^{(i)}_0 &=&
\displaystyle\dfrac{1}{b_l} 
\sum_{
k\equiv l+\omega_{i+1}\atop\mbox{\scriptsize (mod $Q$)}} \prod_{0<j<k} [a_{jk}]
\sum_{\nu =0}^{n-1} 
\dfrac{\mbox{[}v-u +a_{0\nu} \mbox{]}_\omega}{\mbox{[}v-u\mbox{]}} \\
&\times& \displaystyle\oint_{C_\nu} \prod_{j=1}^{n-1} 
\dfrac{dz_j}{2\pi\sqrt{-1}} \prod_{j=0\atop j\neq\nu}^{n-1} 
\dfrac{\mbox{[}a_{0j}\mbox{]}_\omega}{\mbox{[}a_{\nu j}\mbox{]}}
f(v_{j+1}-v_j , 1-a_{\nu j}) A^{(i)}_{l,k}(v; v_1, \cdots , v_{n-1}). 
\end{array}
\label{eq:g_i,0}
\end{equation}
Here, the integral contour $C_\nu$ is chosen such that 
$$
|z_j|=\left\{ \begin{array}{ll} 
x^j (|z|+j\varepsilon ) & (1\leqq j\leqq \nu) \\
x^j (|z|-(n-j)\varepsilon ) & (\nu +1\leqq j\leqq n-1) 
\end{array} \right. 
$$
where $\varepsilon >0$ is a very small positive number. 

Let us denote the r.h.s. of (\ref{eq:g_i,0}) by $H^{(i)}_l$. 
As noted in the previous section, the trace on ${\cal H}^{(i)}_{l,k}$ 
should be taken on ${\cal F}^{(i)}_{\sigma^{-\mu}(l),\sigma^{-\mu}(k)}$ 
for $\mu >0$. Thus, $\langle g \rangle ^{(i)}_\mu$ can be reduced to 
$H^{(i)}_{\sigma^{-\mu}(l)}$. 
Let 
$$
B^{(i)}_{l,k}(v,u):= x^{n\left(\frac{r-1}{r}|k^2|-2\langle k,l\rangle 
+\frac{r}{r-1}|l|^2\right)} x^{2a_{0\,n-1} (v-u)} \tilde{G}_a, ~~~~
\tilde{G}_a=\prod_{j=1}^{n-1} [a_{0j}]_\omega \prod_{0<j<k} [a_{jk}]. 
$$
Consider the following sum 
$$
S^{(i)}(v,u):=\dfrac{[0]_\omega}{[v-u]}\sum_{\mu =0}^{n-1} \sum_{
k\equiv l+\omega_{i+1}\atop\mbox{\scriptsize (mod $Q$)}}
B^{(i)}_{\sigma^{-\mu}(l),\sigma^{-\mu}(k)}(v,u), 
$$
and take the limit $u\rightarrow v$\footnote{
Belavin's $({\mathbb Z}/n{\mathbb Z})$-symmetric model does not have 
the parameter $u$ so that all the physical quantities should be independent 
of $u$. Thus, we set $u\rightarrow v$ here, in order to avoid some 
difficulty. }. Then we have 
\begin{equation}
\lim_{u\rightarrow v} S^{(i)}(v,u)=\omega^{i+1} b_l
\dfrac{(\omega x^{2r}; x^{2r})_\infty 
(\omega^{-1} x^{2r}; x^{2r})_\infty}{(x^{2r}; x^{2r})_\infty^2} 
\dfrac{(x^2; x^2)_\infty (x^{2n}; x^{2n})_\infty^n}{
(\omega ; x^2)_\infty (\omega^{-1}x^2; x^2)_\infty}. 
\label{eq:S-lim}
\end{equation}
This can be confirmed by comparing the series-expansion 
in $x$ of both sides order by order. 

Here we cite the sum formula from \cite{AJMP}: 
\begin{equation}
\sum_{\nu =0}^{n-1} \prod_{j=0\atop j\neq\nu}^{n-1} 
\dfrac{f(v_{j+1}-v_j, 1-\pi_{\nu j})}{[\pi_{\nu j}]}=0. 
\label{eq:sum=0}
\end{equation}
This can be derived by applying the Liouville's second theorem to 
the following elliptic function 
$$
F(w)=\prod_{j=0}^{n-1} \dfrac{[v_{j+1}-v_j-\frac{1}{2}+w-\pi_j ]}
{[v_{j+1}-v_j-\frac{1}{2}][w-\pi_j ]}. 
$$
On eq. (\ref{eq:g_i,0}), 
the contour for $z_{1}$-integral are common for all $\nu$ 
except for $\nu =0$. Thus, by using (\ref{eq:sum=0}), 
$H^{(i)}_l$ can be evaluated by the residue 
at $z_{1}=x^{1+2u}\rightarrow x^{1}z$. The resulting 
$(n-2)$-fold integral has the same structures of both 
the integrand and the contour as the original $(n-1)$-fold one, 
except for the number of integral variables by one. 
Thus, we can repeat this evaluation procedure $n-1$ times to find 
\begin{equation}
H^{(i)}_l\sim \dfrac{1}{[v-u]} \dfrac{1}{b_l} 
\dfrac{B^{(i)}_{l,k}}{(x^{2n}; x^{2n})_\infty^{n-1}}, 
\label{eq:H-lim}
\end{equation}
at $u\sim v$. Substituting (\ref{eq:S-lim}) and (\ref{eq:H-lim}) 
into 
\begin{equation}
\langle g \rangle^{(i)} =\dfrac{1}{\chi^{(i)}} 
\sum_{\mu =0}^{n-1} H^{(i)}_{\sigma^{-\mu}(l)}, 
\end{equation}
we reproduce the expression for 
the spontaneous polarization (\ref{eq:SPn}) originally obtained 
in \cite{SPn}.

\section{Concluding remarks}

In this paper we constructed a free field representation method 
in order to obtain correlation functions of Belavin's 
$({\mathbb Z}/n{\mathbb Z})$-symmetric model. The essential point was 
to find a free field representation of the tail operator $\Lambda_a^{a'}$, 
the nonlocal operator which intertwine 
$({\mathbb Z}/n{\mathbb Z})$-symmetric model and $A^{(1)}_{n-1}$ model. 
As a consistency check, we perform $(n-1)$-fold integrals and traces 
for one-point function, to reproduce the expression of the 
spontaneous polarization originally obtained in \cite{SPn}. 

There are some related works concerning the eight-vertex model 
and its higher spin version. A bootstrap approach 
for the eight-vertex model was presented in \cite{qXYZ}. 
The vertex operators of the eight-vertex model 
with some special values of $r$ were directly bosonized in \cite{Shi4}. 
A free field representation method for form factors of the 
eight-vertex model was constructed in \cite{La}. 
A higher spin generalization of the free field representation method 
was achieved in \cite{KKW}. 
As for $({\mathbb Z}/n{\mathbb Z})$-symmetric model, 
it is important to consider the extension to the form factor problem, 
or the application to the fused model. 
We wish to address these problems in the future.

\section*{Acknowledgements} 

We would like to thank K. Hasegawa, R. Inoue, M. Jimbo, H. Konno, 
M. Lashkevich, T. Miwa, A. Nakayashiki, M. Okado, Ya. Pugai, J. Shiraishi 
and Y. Yamada for discussion and their interest in the present work. 
This work was supported by Grant-in-Aid for Scientific Research 
(C) 19540218 from the Japan Society for the Promotion of Science.

\end{document}